\documentclass[10pt]{article}

\usepackage{amsfonts} 
\usepackage{amsmath,hhline,latexsym,mathrsfs}
\usepackage{amssymb}
\usepackage{relsize}
\DeclareMathAlphabet\mathbfcal{OMS}{cmsy}{b}{n}

\usepackage[latin1]{inputenc}

\usepackage{hyperref}

\usepackage{color}
\usepackage{graphicx}
\usepackage{comment}

\newtheorem{theorem}{\sc Theorem}[section]
\newtheorem{lemma}{\sc Lemma}[section]

\newtheorem{proposition}{\sc Proposition}[section]

\newtheorem{remark}{Remark}

\newcommand{\eps}{\varepsilon}

\newcommand{\N}{\mbox{$I \kern -4pt N$}}
\newcommand{\Q}{\mbox{$Q \kern -8pt I$}}
\newcommand{\R}{\mbox{$I \kern -4pt R$}}
\newcommand{\C}{\mbox{$C \kern -8pt I$}}

\def\R{{\bf R}}
\usepackage[english,ruled,vlined]{algorithm2e}

\title{Asymptotic analysis of an advection-diffusion equation involving interacting boundary and internal layers }

\author{
\textsc{Youcef Amirat}\thanks{Laboratoire de Math\'ematiques Blaise Pascal, Universit\'e Clermont Auvergne, 
UMR CNRS 6620, Campus universitaire des C\'ezeaux, 3, place Vasarely, 63178, Aubi\`ere, France.  E-mail: {\tt youcef.amirat@uca.fr.}}
\quad 
\and
\textsc{Arnaud M\"unch}\thanks{Laboratoire de Math\'ematiques Blaise Pascal, Universit\'e Clermont Auvergne, 
UMR CNRS 6620, Campus universitaire des C\'ezeaux, 3, place Vasarely, 63178, Aubi\`ere, France. E-mail: {\tt arnaud.munch@uca.fr.}}}

\begin{document}

\maketitle	

\begin{abstract}
As $\eps$ goes to zero, the unique solution of the scalar advection-diffusion equation $y^{\eps}_t-\eps y^{\eps}_{xx} + M y^{\eps}_x=0$, $(x,t)\in (0,1)\times (0,T)$ submitted to Dirichlet boundary conditions exhibits a boundary layer of size $\mathcal{O}(\eps)$ and an internal layer of size $\mathcal{O}(\sqrt{\eps})$. If the time $T$ is large enough, these thin layers where the solution $y^{\eps}$ displays rapid variations intersect and interact each other. Using the method of matched asymptotic expansions, we show how we can construct an explicit approximation $\widetilde{P}^\eps$ of the solution $y^\eps$ satisfying $\Vert y^{\eps}-\widetilde{P}^\eps\Vert_{L^\infty(0,T; L^2(0,1))}=\mathcal{O}(\eps^{3/2})$ and $\Vert y^{\eps}-\widetilde{P}^\eps\Vert_{L^2(0,T; H^1(0,1))}=\mathcal{O}(\eps)$, for all $\eps$ small enough. 
 \end{abstract}




\section{Introduction. Problem statement}\label{sec:intro}

Let  $T>0$ and $Q_T:=(0,1)\times (0,T)$.  This work is concerned with the scalar advection-diffusion equation 
\begin{equation}
\label{eq:transport}
	\left\{
		\begin{aligned}
   			& y^{\eps}_{t} - \eps y^{\eps}_{xx}+My^{\eps}_x=0, 				& (x,t)\in Q_T,    \\
   			& y^{\eps}(0,t)  = v^{\eps}(t), \; y^{\eps}(1,t)=0,		  		& t\in  (0,T), \\
   			& y^{\eps}(x, 0) = y^{\eps}_0(x), 						  		& x\in  (0,1),
   		\end{aligned} 
 	\right.
\end{equation}
where $\eps \in (0,1)$ is the diffusion coefficient and $M\in \mathbb{R}^\star$ is the transport coefficient. The function $y^{\eps}_0\in H^{-1}(0,1)$ is the initial data, $v^{\eps}\in L^2(0,T)$ is the boundary data, and $y^{\eps}=y^{\eps}(x,t)$ is the associated state. 

For any $y^{\eps}_0$ in $H^{-1}(0,1)$ and $v^{\eps}$ in $L^2(0,T)$, there exists a unique solution $y^{\eps}$ to (\ref{eq:transport}), with the regularity $y^{\eps}\in L^2(Q_T)\cap C([0,T];H^{-1}(0,1))$.

We are interested in this work with a precise asymptotic description of the solution $y^{\eps}$ when $\eps$ is small. As a first motivation, we mention that the system (\ref{eq:transport}) can be seen as a simple example of complex models where the diffusion coefficient is small compared to the others. Actually, as discussed in \cite{temam2000}, the model problem \eqref{eq:transport} is an embedded system of the Navier-Stokes system with non-characteristic boundary condition and viscosity coefficient equals to $\eps$. A second motivation comes from the numerical approximation of (\ref{eq:transport}) that may be not straightforward for small values of $\eps$ (we refer to \cite{eymard},\cite{stynes1989}). A third motivation comes from the asymptotic controllability property of (\ref{eq:transport}) studied in \cite{CoronGuerrero2005} and which exhibits surprising behaviors, leaving many open questions. 

\begin{figure}[!http]
\begin{center}
\includegraphics[scale=0.65]{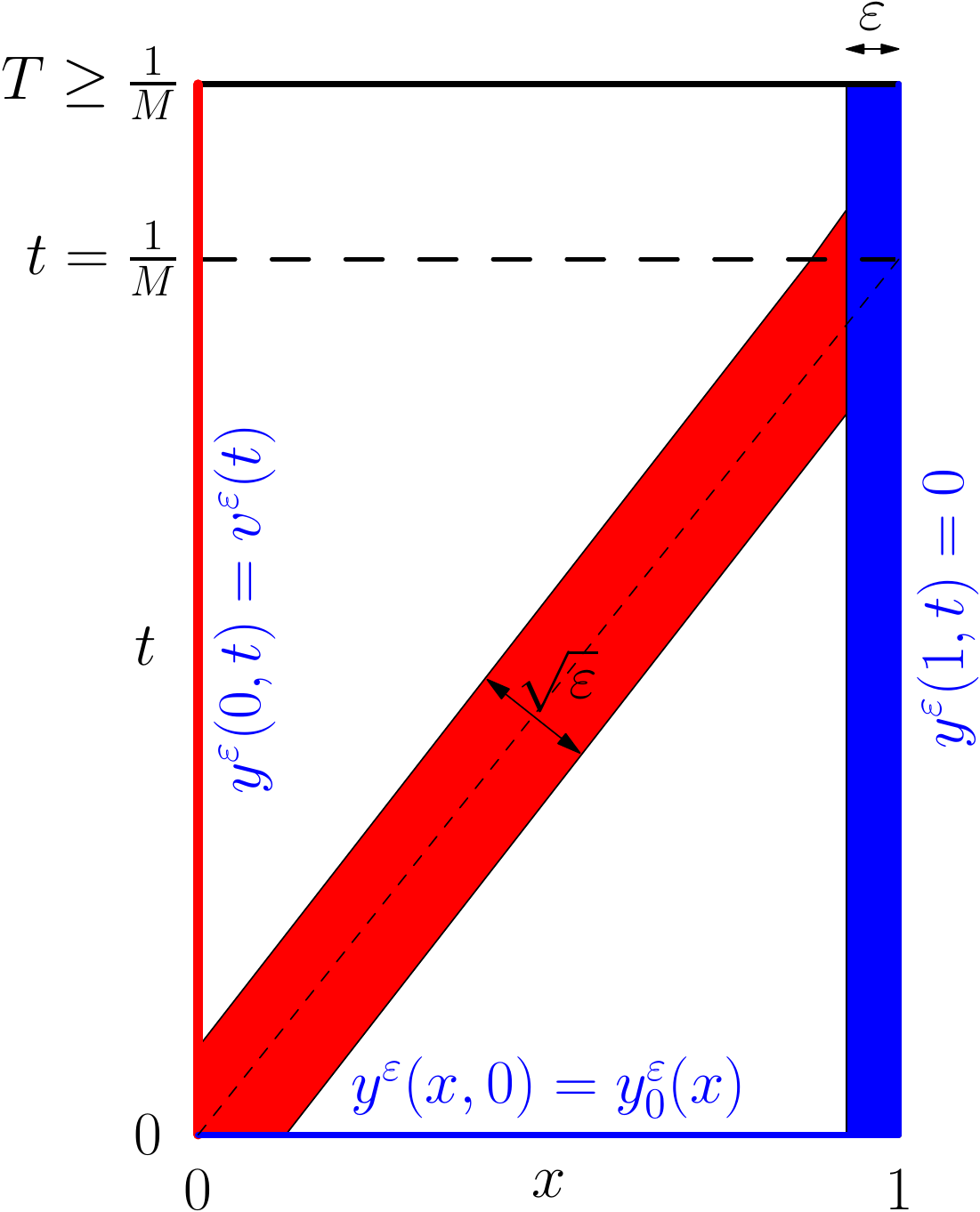}
\caption{Internal (red) and boundary (blue) layer zones for $y^{\eps}$  in the case $M>0$.}\label{boundarylayer}
\end{center}
\end{figure}

The limit of the solution $y^{\eps}$ has been studied in \cite{CoronGuerrero2005}, assuming that the initial condition does not depend on $\eps$, that is $y_0^{\eps}=y_0$. Precisely, it is shown that, if $(v^\eps)_{(\eps>0)}$ is a sequence of functions in $L^2(0,T)$ such that, for some $v\in L^2(0,T)$, $v^\eps\rightharpoonup v$ in $L^2(0,T)$ weak, as $\eps\to 0^+$, then $y^{\eps}\rightharpoonup y$ in $L^2(Q_T)$ weak, as $\eps\to 0^{+}$, where $y\in C([0,T];L^2(0,1))$ is the weak solution of the following transport equation
\begin{equation*}
\label{eq:transportlimit}
\left\{
\begin{aligned}
	&  y_{t} +My_x=0, 			& 	(x,t)\in Q_T,   \\
	&  y(0,t)  = v(t),     \quad  \textrm{if}\quad M>0  &   t\in (0,T),   \\
	&  y(L,t)  = 0,   \quad  \textrm{if}\quad M<0&   t\in (0,T),  \\
	&  y(x, 0) = y_0(x),  & x\in (0,1).
\end{aligned}
\right.
\end{equation*}
Despite its apparent simplicity, the asymptotic analysis of system (\ref{eq:transport}) with respect to $\eps$ is not straightforward. Take for instance $M>0$. The reason is that, as $\eps$ goes to zero, the solution $y^{\eps}$ exhibits a boundary layer of size $\mathcal{O}(\eps)$ at $x=1$ (blue part on Figure \ref{boundarylayer}) but also an internal layer (also called shock layer in some cases) of size $\mathcal{O}(\sqrt{\eps})$ along the characteristic $\{(x,t)\in Q_T, x-M t=0\}$ (red part on Figure \ref{boundarylayer}).
Thus, two distinct layers, with different sizes, appear and interact in the neighborhood of the point $(x,t)=(1,1/M)$, assuming $T\geq 1/M$.
  
The boundary layer at the boundary $x=1$ occurs as soon as the initial condition $y_0^{\eps}$ is different from zero. On the other hand, whatever be the regularity of the initial condition $y_0^{\eps}$ and Dirichlet condition $v^\eps$, the internal layer along the characteristic occurs if $y_0^{\eps}$ and  $v^{\eps}$ do not satisfy appropriate compatibility conditions at the point $(0,0)$. For instance, if both $v^{\eps}=v$ and $y_0^\eps=y_0$ are independent of $\eps$, these compatibility conditions read as follows: 
\begin{equation}\label{compability_condition}
M^p y_0^{(p)}(0)+(-1)^{p+1}v^{(p)}(0)=0, \quad \forall p\in \mathbb{N}.
\end{equation}
We refer to \cite{amirat_munch}. Assuming such conditions, the asymptotic analysis of (\ref{eq:transport}) has been done in \cite{amirat_munch}. More precisely, assuming that the initial condition is independent of $\eps$ and that the function $v^{\eps}$ is given in the form $v^{\eps}=\sum_{k=0}^m \eps^k v^k$,  an asymptotic approximation $w_m^{\eps}$ of the solution $y^{\eps}$ is constructed in \cite{amirat_munch}. The method of matched asymptotic expansions is used to define an outer solution (out of the boundary layer) and an inner solution. Upon regularity assumptions on the functions $v^k$, $k=0,\dots,m$ and $y_0^{\eps}=y_0$, $w_m^{\eps}$ is shown to be a regular and strong convergent approximation of $y^{\eps}$, as $\eps\to 0^+$. For any $m\in \mathbb{N}$, the error estimate is as follows
\begin{equation*}
\label{i1}
\Vert y^\varepsilon(\cdot,t) - w_m^\varepsilon(\cdot,t) \Vert_{L^2(0,1)} \leq c_m \varepsilon^{\frac{2m+1}{2}\gamma} +
c_m \left(\varepsilon^{\frac{1}{2}} + \varepsilon^{\frac{(2m+3)\gamma}{2}}\right)
e^{-\frac{M^2}{2 \eps^\gamma} t},
\quad \forall t\in [0,T],
\end{equation*}
for some constant $c_m$ independent of $\eps$ and $\gamma\in (0,1/2]$. The function $w_m^{\eps}$, sum of solutions of transport equations, explicit, can therefore be used for numerical purposes. The estimate involves the initial layer corrector, exponentially small with respect to $\eps$ for $t>0$. Moreover, assuming that $y_0^{\eps}$ is a Gevrey function of order $1/2$ in $[0,1]$, and that the $v^k$ functions are polynomials, the constant $c_m$ is uniformly bounded with respect to $m$ allowing to pass to the limit, as $m\to \infty$, 
with $\eps$ small enough but fixed. This leads to the following decomposition 
\begin{equation*}
\label{i2}
y^\eps(x,t)=w^{\eps}(x,t)+ \theta^\eps(x,t), \quad (x,t)\in Q_T,
\end{equation*}
where $w^{\eps}$ is an infinite sum of explicit solutions of transport equations, and $\theta^\eps$ is the initial layer corrector, defined as the solution of a non-homogeneous advection-diffusion equation of the form
(\ref{eq:transport})
satisfying $\Vert \theta^{\varepsilon}(\cdot,t)\Vert_{L^2(0,1)}  \leq c \varepsilon^{\frac{1}{2}}
e^{-\frac{M^2}{2 \eps^\gamma} t}$, for all $t\in [0,T]$, for some constant $c$ independent of $\eps$.
 
The main purpose of this work is to reproduce (partially) the analysis done in \cite{amirat_munch}, relaxing the compatibility conditions \eqref{compability_condition}. These conditions are notably very restrictive for the underlying null controllability problem where $v$ is sought such that $y^{\eps}=y^{\eps}(y_0,v)$ vanishes at any controllability time $T$ (see \cite{CoronGuerrero2005},\cite{amirat_munch_sinica}). The violation of the compatibility conditions create a thin inner region (called internal layer) in the vicinity of the characteristic $\{(x,t)\in Q_T, x-M t=0\}$ where the solution $y^{\eps}$ exhibits rapid variations. The internal layer will intersect the boundary layer living along $x=1$, when $T\geq 1/M$. This requires to incorporate to our analysis the internal layer appearing along the characteristic and notably, to discuss how this internal layer interacts with the boundary layer at $x=1$.

Interaction of shock layers and boundary layers is a well-know phenomenon, for instance in fluid mechanics. We refer to \cite{adamson1980} and the references therein. Asymptotic analysis with respect to a small parameter of a boundary value problem involving such interaction is however quite challenging. To our knowledge, few analysis have been performed, mainly on simple cases. Influence of mutual interaction of layers is notably discussed, based on heuristic scaling arguments, in \cite{howes1989} (see also \cite{Shao2007}). The first example considered in \cite{howes1989} is  
\begin{equation}
\label{Pb_elliptic}
\left\{
\begin{aligned}
& u^{\eps}_x(x,y)-\eps \Delta u^{\eps}(x,y)=0, & (x,y)\in (0,1)\times (-1,1), \\
& u^{\eps}(0,y)=f(y), \quad u^{\eps}(x,-1)=u^{\eps}(x,1)=u^{\eps}(1,y)=0, & x\in [0,1], y\in [-1,1],
\end{aligned}
\right.
\end{equation}
where $f:[-1,1]\to \mathbb{R}$ is a piecewise constant, discontinuous at $y=0$.  This develops a shock layer around the line $\{y=0\}$, intersecting at the point $(1,0)$, the boundary layer living along the orthogonal line $\{x=1\}$. This example is inspired from \cite{eckhaus66} where a rigorous constructive asymptotic analysis is given and leads to an approximation $w^{\eps}$ satisfying the uniform property $\vert u^{\eps}(x,y)-w^{\eps}(x,y)\vert=\mathcal{O}(\sqrt{\eps})$ for all $(x,y)\in [0,1]\times [-1,1]$. This estimate is obtained using a maximum principle. Later on, using similar technics, the asymptotic analysis of the system (\ref{eq:transport}) has been discussed in \cite{Bobisud67}; precisely, assuming $v\in C^3([0,T])$ and $y_0\in C^4([0,1])$, a continuous explicit function converging to $y^{\eps}$ for the uniform norm with a rate $\sqrt{\eps}$ is determined. 

It is also worth mentioning the works \cite{Shih1996,Shih2001} where the asymptotic analysis of the system (\ref{eq:transport}) defined over $\mathbb{R}^+\times (0,T)$ is performed, for $M>0$. In that case, there is not more boundary layer. Using the matching asymptotic method, it is shown that the solution in the internal layer can be represented by iterated integral of the complementary error function erfc. In particular, the analysis provides the exact solution for the integral representation of the solution $y^{\eps}$ when the initial condition $y_0$ and boundary data $v$ are polynomials. In the general case, estimate for the uniform norm are obtained using maximum principles. For instance, the shock layer (appearing when $v(0)\neq y_0(0)$) is analyzed in \cite{Shih1994} and leads to an error for the uniform norm of the order $\eps^{1/2}$; the angular layer (appearing when $v(0)=y_0(0)$ but $v^{(1)}(0)\neq -M y_0^{(1)}(0)$) leads to an error of the order $\eps^{3/2}$. 

Assuming the data $v$ and $y_0$ in $C^4(0,T)$ and $C^4(0,1)$ respectively, we construct in this work an explicit $C^1$-approximation of $y_{\eps}$ leading to error estimate in $L^\infty(0,T,L^2(0,1))$ and $L^2(0,T;H_0^1(0,1))$. We use the method of matched asymptotic expansions together with energy estimates. The analysis combines in an appropriate way the description of the internal layer given in \cite{Shih1994, Shih1996} to the description of the boundary layer given in \cite{amirat_munch}.  The document is organized as follows. In Section \ref{section2}, we employ the method of matched asymptotic expansions to construct a linear combination of three distinct expansions: a first outer expansion, defined as the sum of the functions $y^k$, $k=0,\dots,3$, solution of pure advection equations, aims to approximate the solution $y^\eps$, far away from the boundary and internal layers. A second inner expansion, defined as the sum of the function $W^k$, $k=0,\dots,3$, solution of pure diffusion equations, aims to approximate the solution $y^{\eps}$ in a neighborhood of size $\mathcal{O}(\sqrt{\eps})$ of the first characteristic. Eventually, a third inner expansion, defined as the sum of the function $Y^k$, $k=0,\dots,3$, solution of ordinary differential equations, aims to approximate the solution $y^{\eps}$ in a neighborhood of size $\mathcal{O}(\eps)$ of the boundary $x=1$. A composite technic which consists in adding these three expansions and then subtracting their commons parts leads to a first sequence $(P^\eps)_{(\eps>0)}$. As $\eps$ goes to zero, it turns out that this sequence exhibits an artificial corner layer in the neighborhood of the point $(x,t)=(0,0)$ leading to an unsatisfactory approximation result, namely $\Vert P^\eps-y^{\eps}\Vert_{L^\infty(0,T,L^2(0,1))}=\mathcal{O}(\sqrt{\eps})$. For this reason, in Section \ref{section3}, using precise descriptions of the internal layer given in \cite{Shih1994, Shih1996}, we slightly modify the sequence $(P^\eps)_{(\eps>0)}$ and obtain a second sequence $(\widetilde{P}^\eps)_{(\eps>0)}$.  Then, in Section \ref{section_convergence}, using energy estimates for the advection-diffusion equation (\ref{eq:transport}) together with precise $L^1(L^2)$ estimates of some remainder terms, we prove the convergence of $(\widetilde{P}^{\eps}-y^{\eps})_{(\eps>0)}$ to zero for the $L^\infty(L^2)$-norm with the rate $\eps^{3/2}$, see Theorem \ref{convergence}. Assuming moreover that $y_0(1)=y_0^{(1)}(1)=0$, we also show the convergence of the sequence $(\widetilde{P_x}^{\eps}-y^{\eps}_x)_{(\eps>0)}$  to zero for the $L^2(Q_T)$-norm with the rate $\eps$, see Theorem \ref{convergenceL2H1}. Then, as an application of the $L^\infty(L^2)$-estimate, we show that the $L^2(0,1)$-norm of the solution $y^{\eps}(\cdot,1/M)$  at time $1/M$ associated to $v=0$ decays polynomially with respect to $\eps$ and depends only of the value of the initial condition at $x=0$, see Proposition \ref{polynomial_decay}. Section \ref{conclusions} concludes with some perspectives. 

As far as we know, this study is the first detailed  asymptotic analysis of a boundary value problem involving two interacting singular layers.

\section{Matched asymptotic expansions and approximate solutions}
\label{section2}

Let us consider the problem 
\begin{equation}
\label{eq:2.1}
\left\{
\begin{aligned}
& y^{\varepsilon}_t - \varepsilon y^{\varepsilon}_{xx} + M y^{\varepsilon}_x=0, & (x,t) \in Q_T,\\
& y^{\varepsilon}(0,t)=v(t),  \quad y^{\varepsilon}(1,t)=0, & t\in (0,T),\\
& y^{\varepsilon}(x,0)=y_0(x), & x\in (0,1),
\end{aligned}
\right.
\end{equation}
where $y_0$ and $v$ are given functions, $M>0$, and $T\geq 1/M$ . We construct an asymptotic approximation of the solution $y^{\varepsilon}$ of~\eqref{eq:2.1} by using the method of matched asymptotic expansions. We refer to \cite{ColeBook81, SanchezPalencia, VanDyke, eckhausbook79} for a general presentation of the method. We assume that the initial condition $y^{\varepsilon}(x,0)$ and the boundary condition $y^{\varepsilon}(0,t)$ are independent of $\eps$. In view of the linearity of (\ref{eq:2.1}), the procedure is very similar for $y^{\varepsilon}(\cdot,0)$ and $y^{\varepsilon}(0,\cdot)$ of the form, respectively, $y^{\varepsilon}(\cdot,0)=\sum_{k=0}^m \varepsilon^{k} y_0^{k}$ and
$y^{\varepsilon}(0,\cdot)=\sum_{k=0}^m \varepsilon^{k} v^{k}$. Note that the case $M<0$ can be treated similarly. 



As mentioned in the introduction, the solution $y^{\eps}$ exhibits two inner regions: an internal layer located along the characteristic $\{(x,t)\in Q_T, x-Mt=0\}$ and a boundary layer living along $x=1$. The internal layer is of size $\mathcal{O}(\varepsilon^{1/2})$ while the boundary layer is of size $\mathcal{O}(\varepsilon)$. The outer region is the subset of $(0,1)$ consisting of the points far from the internal and the boundary layers, it is of $\mathcal{O}(1)$ size. The occurrence of these three distinct regions require to introduce three distinct asymptotic expansions. 
The first one, the so-called outer expansion, lives far away from the inner regions and is given by   
\begin{equation*}
\sum_{k=0}^m \varepsilon^{k} y^{k}(x,t), \quad (x,t)\in Q_T,   \quad x-Mt\neq 0, \quad x<1,
\end{equation*}
for some $m\in \mathbb{N}^\star$. A second one, the so-called first inner expansion, living in the neighborhood of $\{(x,t)\in Q_T, x-Mt=0\}$ is given by 
\begin{equation*}
\sum_{k=0}^{m} \eps^\frac{k}{2} W^{k/2}(w,t), \quad w=\frac{x-Mt}{\eps^{1/2}} \in \biggl(-\frac{Mt}{\eps^{1/2}},\frac{1-Mt}{\eps^{1/2}}\biggr), \quad t\in(0,T).
\end{equation*}
Last, a third one, the so-called second inner expansion, living along $x=1$, is given by 
\begin{equation*}
\sum_{k=0}^{m} \varepsilon^{k/2} Y^{k/2}(z,\tau,t), \quad z=\frac{1-x}{\varepsilon} \in (0,\varepsilon^{-1}), \quad \tau=\frac{\frac{1}{M}-t}{\sqrt{\eps}}.
\end{equation*}
In particular, these expansions make appear several variables, at different scales, namely, 
$$
x,\quad t, \quad z=\frac{1-x}{\eps}, \quad w=\frac{x-Mt}{\sqrt{\eps}}, \quad \tau=\frac{\frac{1}{M}-t}{\sqrt{\eps}}.
$$
We will construct outer and inner expansions which will be valid in the so-called outer and inner regions, respectively. There are intermediate regions between the outer region and the inner regions, with size $\mathcal{O}(\varepsilon^\gamma)$, $\gamma \in (0,1)$. To construct an approximate solution we require that inner and outer expansions coincide in each intermediate region, then some conditions must be satisfied in that region by the corresponding inner and outer expansions. These conditions are the so-called matching asymptotic conditions.

The strategy is as follows. We first identify the functions $y^k$, $k=0,\dots,m$ in the outer region. Then, we identify the functions $W^{k/2}$, $k=0,\dots,m$ of the first inner expansion satisfying the matching conditions (with the $y^k$). 
This allows to define an expansion, denoted by $p^{\eps}$, valid far away from $x=1$, as a linear combination of the functions $y^k$ and $W^{k/2}$. Then, we identify the functions $Y^k$, $k=0,\dots,m$ of the second inner expansion satisfying the matching conditions with $p^\eps$. Eventually, we define an expansion, denoted by $P^{\eps}$, valid in the whole domain, as linear combination of the function $Y^k$ and the function $p_\eps$, and supposed to be an approximation of $y^{\eps}$. In this work, we shall take $m=3$.

\subsection{Outer expansion}

Putting $\displaystyle \sum_{k=0}^m \varepsilon^{k} y^{k}(x,t)$ into equation~$\eqref{eq:2.1}_1$, the identification of the powers of $\varepsilon$ yields
\begin{align*}
& \varepsilon^0: \quad y_t^0 +M y_{x}^0=0, \\
& \varepsilon^{k}: \quad y_t^{k} +M y_{x}^{k}=y^{k-1}_{xx}, \quad \mbox{ for any } 1\leq k\leq m. 
\end{align*}
Taking the initial and boundary conditions into account we define $y^0$ and $y^{k}$ $(1\leq k\leq m)$ as functions  satisfying the transport equations, respectively,
\begin{equation}
\label{eq2.2:y0}
\left\{
\begin{aligned}
& y_t^0 +M y_{x}^0=0, & (x,t)\in Q_T, \\
& y^0(0,t)=v(t), & t\in (0,T), \\
& y^0(x,0)=y_0(x), & x\in (0,1), 
\end{aligned}
\right. 
\quad\qquad \left\{
\begin{aligned}
& y_t^{k} +M y_{x}^{k}=y^{k-1}_{xx}, & (x,t)\in Q_T, \\
& y^{k}(0,t)=0, & t\in (0,T), \\
& y^{k}(x,0)=0, & x\in (0,1). 
\end{aligned}
\right.
\end{equation}
The solution $y^0$ is given by
\begin{equation}
\label{2.y0def}
y^0(x,t)=
\left\{
\begin{aligned}
& y_0(x-Mt), & x> Mt, \\
& v\left(t-\frac{x}{M}\right), & x < Mt. 
\end{aligned}
\right. 
\end{equation}
Using the method of characteristics we find that, for any $1\leq k\leq m$, the solution $y^k$ is given by
\begin{equation*}
\label{2.ykdef}
y^{k}(x,t)=
\left\{
\begin{aligned}
&  \int_0^t   y^{k-1}_{xx}(x+(s-t)M,s) ds,    & x > Mt, \\
& \int_0^{x/M}   y^{k-1}_{xx}(s M, t-\frac{x}{M}+s)ds, &   x < Mt.
\end{aligned}
\right.
\end{equation*}
\begin{remark}
\label{rem2.1}
We may determine explicitly the function $y^k$. For instance, we verify that we have 
\begin{equation*}
\label{2.2y1def}
y^1(x,t)=
\left\{
\begin{aligned}
& t \,y^{(2)}_0(x-M t),    & x > Mt, \\
& \frac{x}{M^3} v^{(2)}\left(t-\frac{x}{M}\right), &   x < Mt,
\end{aligned}
\right.
\end{equation*}
and
\begin{equation*}
\label{2.2y2def}
y^2(x,t)=
\left\{
\begin{aligned}
& \frac{t^2}{2} \,y^{(4)}_0(x-M t),    & x > Mt, \\
&-\frac{2x}{M^5} v^{(3)}\left(t-\frac{x}{M}\right)  + \frac{x^2}{2M^6} v^{(4)}\left(t-\frac{x}{M}\right), &   x < Mt.
\end{aligned}
\right.
\end{equation*}
Here and in the sequel, $f^{(i)}$ denotes the derivative of order $i\in \mathbb{N}$ of the real function $f$. $\hfill\Box$
\end{remark}

\subsection{Inner expansion along the characteristic}
\label{innerexp}

Now we turn back to the construction of the first inner expansion. In the sequel, we shall use the error function $erf:\mathbb{R}\to [-1,1]$ defined as $erf(y)=\frac{2}{\pi}\int_0^y e^{-s^2}ds$. 
It satisfies notably the estimates 
$$
\sqrt{1-e^{-y^2}}\leq erf(y)\leq \sqrt{1-e^{-\frac{4y^2}{\pi}}}, \quad \forall y\geq 0
$$
proved in \cite{Chu} so that $erf(y)\to 1$ as $y\to \infty$. In particular, these estimates imply
\begin{equation}
\frac{1}{2}e^{-\frac{4y^2}{\pi}}\leq 1-erf(y)\leq e^{-y^2}, \quad \forall y\geq 0. \label{estimate_chu}
\end{equation}
We will also use in the sequel the asymptotic behavior of the erf function for large $y>0$:
\begin{equation}
erf(y)=1-e^{-y^2} \biggl(\frac{1}{\sqrt{\pi}y }+\mathcal{O}\biggl(\frac{1}{y^3}\biggr) \biggr).  \label{behavior_erf}
\end{equation}
Since $erf(-y)=-erf(y)$ for all $y\in \mathbb{R}$, similar relations holds for $y<0$, in particular $erf(y)\to -1$ as $y\to -\infty$. Eventually, we also introduce the complementary error function  
$erfc:\mathbb{R}\to [0,2]$ defined as $erfc(y)=1-erf(y)$.

Putting $\displaystyle \sum_{k=0}^{m} \varepsilon^\frac{k}{2} W^{k/2}(w,t)$ into equation~$\eqref{eq:2.1}_1$, the identification of the powers of $\varepsilon$ yields
\begin{align*}
W_t^{k/2}(w,t) - W^{k/2}_{ww}(w,t)=0,   \quad \mbox{ for any } 0\leq k\leq m. 
\end{align*}
so that $W^k$ solves, for all $0\leq k\leq m$, the heat equation. To get the asymptotic matching conditions we write that, for any fixed $t$ and large $w$,
\begin{align*}
& W^0(w,t) + \eps^{1/2} W^{1/2}(w,t) + \varepsilon W^1(w,t) + \cdots + \varepsilon^{m/2} W^{m/2}(w,t) \\
& =y^0(x,t) +   \varepsilon y^1(x,t) + \cdots + \varepsilon^{m} y^{m}(x,t) + \mathcal{O}(\varepsilon^{(m+1)}).
\end{align*}
Rewriting the right-hand side of the above equality in terms of $w$, $t$ and using Taylor expansions we have
\begin{align*}
& W^0(w,t) + \sqrt{\varepsilon} W^{1/2}(w,t) + \varepsilon W^1(w,t) + \cdots + \varepsilon^{m/2} W^{m/2}(w,t) \\
& = y^0(\sqrt{\eps}w+Mt,t)  + \varepsilon y^1(\sqrt{\eps}w+Mt,t) + \cdots + \varepsilon^{m} y^{m}(\sqrt{\eps}w+Mt,t) + \mathcal{O}(\varepsilon^{(m+1)}) \\
& = y^0(Mt,t) + \sqrt{\eps}wy^0_x(Mt,t)  + \frac{\eps w^2}{2}y^0_{xx}(1,t) +\cdots   
\end{align*}
Therefore, at the first orders, the matching conditions read
\begin{equation}
\label{matching_condition_w}
\begin{aligned}
& W^0(w,t) \sim y^0((Mt)^{\pm},t), \quad \mbox{ as } w\to \pm\infty, \\
& W^{1/2}(w,t) \sim w y^0_x((Mt)^{\pm},t)   , \quad \mbox{ as } w\to \pm\infty, \\
& W^1(w,t) \sim \frac{w^2}{2} y^0_{xx}((Mt)^\pm,t) + y^1( (Mt)^\pm,t), \quad \mbox{ as } w\to \pm\infty,   \\
& W^{3/2}(w,t) \sim \frac{w^3}{3 !} y^0_{xxx}((Mt)^\pm,t) +  w y_x^1( (Mt)^\pm,t), \quad \mbox{ as } w\to \pm\infty.   
\end{aligned}
\end{equation}
%

%
\par\noindent
$\bullet$ We define $W^0$ as a solution of
\begin{equation}
\label{eq2.1:W0}
\left\{
\begin{aligned}
& W^0_t(w,t) - W_{ww}^0(w,t)=0, & (w,t)\in \mathbb{R}\times(0,T), \\
& \lim_{w\to \pm\infty} W^0(w,t)=\lim_{x\to (Mt)^\pm} y^0(x,t), & t\in (0,T).
\end{aligned}
\right.
\end{equation}
In view of (\ref{2.y0def}), we have $\lim_{x\to (Mt)^+} y^0(x,t)=y_0(0)$ and $\lim_{x\to (Mt)^-} y^0(x,t)=v(0)$ so that (\ref{eq2.1:W0}) rewrites
\begin{equation*}
\label{eq2.1:W0p}
\left\{
\begin{aligned}
& W^0_t(w,t) - W_{ww}^0(w,t)=0, & (w,t)\in \mathbb{R}\times(0,T), \\
& \lim_{w\to +\infty} W^0(w,t)=y_0(0), \quad \lim_{w\to -\infty} W^0(w,t)=v(0), & t\in (0,T). 
\end{aligned}
\right.
\end{equation*}
The solution of (\ref{eq2.1:W0}) is not unique. Actually, using the fundamental solution of the heat equation, 
the general form of $W^0$ is as follows: 
\begin{equation*}
W^0(w,t)=(H(\cdot,t)\star g_0)(w)=\frac{1}{\sqrt{4\pi t}} \int_{\mathbb{R}} e^{-\frac{(w-s)^2}{4t}} g_0(s)ds, \quad g_0(w):=\lim_{t\to 0^+} W^0(w,t)
\end{equation*}
with $H:\mathbb{R}\times \mathbb{R}^{+}\to \mathbb{R}^+$ defined by $H(w,t)=(4\pi t)^{-1/2}e^{-w^2/(2t)}$. From this expression, we check that any function $g_0$ satisfying $lim_{w\to -\infty} g_0(w)=v(0)$ and $lim_{w\to \infty} g_0(w)=y_0(0)$ leads to a function $W^0$ fulfilling the prescribed asymptotic behavior in (\ref{eq2.1:W0}). The simplest choice for the initial condition $g_0$ is given by 
\begin{equation}
\label{def.g0}
g_0(w)=
\left\{
\begin{aligned}
&y_0(0), & w\geq 0, \\
& v(0), &  w<0, 
 \end{aligned}
 \right.
 \end{equation}
 leading to the following explicit expression: 
\begin{equation}
\label{def.W0}
\left\{
\begin{aligned}
& W^0(w,t)=\frac{c^+-c^-}{2} erf\biggl(\frac{w}{2\sqrt{t}}\biggr) + \frac{c^++c^-}{2}=c^++ \frac{c^--c^+}{2}erfc\biggl(\frac{w}{2\sqrt{t}}\biggr),\\
& c^+=  y_0(0), \quad c^-=v(0).
\end{aligned}
\right.
\end{equation}
Using the asymptotic behavior of the error function, we easily verify that $W^0$ satisfies the prescribed asymptotic behavior. 
Remark also that $lim_{t\to 0^+} W^0(w,t)$ equals $y_0(0)$ if $w\geq 0$ and $v^0(0)$ if $w<0$. 

\par\noindent
$\bullet$ We define $W^{1/2}$ as a solution of
\begin{equation}
\label{eq2.1:W1d2}
\left\{
\begin{aligned}
& W^{1/2}_t(w,t) - W_{ww}^{1/2}(w,t)=0, & (w,t)\in \mathbb{R}\times(0,T), \\
& \lim_{w\to \pm\infty} \biggl( W^{1/2}(w,t) -  y^0_x((Mt)^\pm,t) \,w\biggr) =0, & t\in (0,T).
\end{aligned}
\right.
\end{equation}
In view of (\ref{2.y0def}), $y_x^0((Mt)^+,t)=(y_0)^{(1)}(0)$ and $y^0_x((Mt)^-,t)=-\frac{1}{M}v^{(1)}(0)$. Then (\ref{eq2.1:W1d2}) rewrites 
\begin{equation*}
\left\{
\begin{aligned}
& W^{1/2}_t(w,t) - W_{ww}^{1/2}(w,t)=0, & (w,t)\in \mathbb{R}\times(0,T), \\
& \lim_{w\to +\infty} \biggl( W^{1/2}(w,t) -  (y_0)^{(1)}(0) w  \biggr) =0, & t\in (0,T),\\
& \lim_{w\to -\infty} \biggl( W^{1/2}(w,t) + \frac{1}{M} v^{(1)}(0) \,w\biggr) =0, & t\in (0,T).
\end{aligned}
\right.
\end{equation*}
Proceeding as before, a solution is given by $W^{1/2}(w,t)=(H(\cdot,t)\star g_{1/2})(w)$ with $g_{1/2}(s)= y_0^{(1)}(0) s$ if $s\geq 0$ and $g_{1/2}(s)=-\frac{1}{M} v^{(1)}(0) \,s$ if $s<0$. Explicitly, we obtain
\begin{equation}
\label{def.W12}
\left\{
\begin{aligned}
& W^{1/2}(w,t) =w\biggl(\frac{d^+-d^-}{2} erf\biggl(\frac{w}{2\sqrt{t}}\biggr)+ \frac{d^++d^-}{2}\biggl) +  (d^+-d^-)\frac{\sqrt{t}}{\sqrt{\pi}} e^{-\frac{w^2}{4t}},\\
& d^+=  y_0^{(1)}(0), \quad d^-=-\frac{1}{M}v^{(1)}(0),
\end{aligned}
\right.
\end{equation}
and check that $W^{1/2}$ satisfies the prescribed asymptotic property.
Remark that $lim_{t\to 0^+} W^{1/2}(w,t)$ equals $y_0^{(1)}(0)w$ if $w>0$ and $-\frac{v^{(1)}(0)}{M} w$ if $w<0$. 
Moreover, 
\begin{equation*}
\lim_{x\to 0^+} W^{1/2}\biggl(\frac{x}{\sqrt{\eps}},0\biggr)-d^+ \frac{x}{\sqrt{\eps}}=0,
\end{equation*}
and
\begin{equation}
\lim_{t\to 0^+} \biggl[W^{1/2}\biggl(\frac{-Mt}{\sqrt{\eps}},t\biggr)-\frac{d^++d^-}{2}\biggl(-\frac{Mt}{\sqrt{\eps}}\biggr)\biggl]=0. \label{limit_t0W12}
\end{equation}


\par\noindent
$\bullet$ We define $W^1$ as a solution of
\begin{equation}
\label{eq2.1:W1}
\left\{
\begin{aligned}
& W^1_t(w,t) - W_{ww}^1(w,t)=0, & (w,t)\in \mathbb{R}\times(0,T), \\
&  \lim_{w\to \pm\infty} \biggl( W^1(w,t) -  \biggl[\frac{w^2}{2} y^0_{xx}((Mt)^\pm,t) + y^1( (Mt)^\pm,t)\biggr]\biggr) =0, & t\in (0,T).
\end{aligned}
\right.
\end{equation}
We have $y^1((Mt)^+,t)=t (y_0)^{(2)}(0)$ and $y^1((Mt)^-,t)=\frac{t}{M^2}v^{(2)}(0)$ for all $t$.  Similarly, 
$y_{xx}^0( (Mt)^+,t )=(y_0)^{(2)}(0)$ and $y_{xx}^0( (Mt)^-,t )=\frac{1}{M^2}v^{(2)}(0)$. Then \eqref{eq2.1:W1} rewrites
\begin{equation*}
\left\{
\begin{aligned}
& W^1_t(w,t) - W_{ww}^1(w,t)=0, & (w,t)\in \mathbb{R}\times(0,T), \\
& \lim_{w\to +\infty} \biggl( W^1(w,t) -  \frac{w^2}{2} y_0^{(2)}(0) \biggr) =t y_0^{(2)}(0), & t\in (0,T),\\
&  \lim_{w\to -\infty} \biggl( W^1(w,t) -  \frac{w^2}{2} \frac{1}{M^2}v^{(2)}(0) \biggr) =\frac{t}{M^2}v^{(2)}(0), & t\in (0,T).
\end{aligned}
\right.
\end{equation*}
A solution is given by $W^1(w,t)=(H(\cdot,t)\star g_1)(w)$ with $g_1(s)= \frac{s^2}{2} y_0^{(2)}(0)$ if $s\geq 0$ and $g_1(s)=\frac{s^2}{2} \frac{1}{M^2}v^{(2)}(0)$ if $s<0$. Explicitly, we obtain
\begin{equation}
\label{def.W1}
\left\{
\begin{aligned}
& W^1(w,t)= \biggl(\frac{w^2}{2}+t\biggr)\biggl(\frac{e^+-e^-}{2} erf\biggl(\frac{w}{2\sqrt{t}}\biggr) +  \frac{e^++e^-}{2}\biggr) + \frac{e^+-e^-}{2} w\sqrt{\frac{t}{\pi}} e^{-\frac{w^2}{4t}},\\
& e^+=y_0^{(2)}(0), \quad e^-=\frac{v^{(2)}(0)}{M^2}.
\end{aligned}
\right.
\end{equation} 

\par\noindent
$\bullet$ We define $W^{3/2}$ as a solution of
\begin{equation}
\label{eq2.1:W3d2}
\left\{
\begin{aligned}
& W^{3/2}_t(w,t) - W_{ww}^{3/2}(w,t)=0, & (w,t)\in \mathbb{R}\times(0,T), \\
&  \lim_{w\to \pm\infty} \biggl( W^{3/2}(w,t) -  \biggl[\frac{w^3}{3 !} y^0_{xxx}((Mt)^\pm,t) +  w y_x^1( (Mt)^\pm,t)\biggr]\biggr)=0, & t\in (0,T).
\end{aligned}
\right.
\end{equation}
We have $y^0_{xxx}((Mt)^+,t)=y_0^{(3)}(0)$, $y_x^1( (Mt)^+,t)=t y_0^{(3)}(0)$, $y^0_{xxx}((Mt)^-,t)=-\frac{1}{M^3}v^{(3)}(0)$ and $y_x^1( (Mt)^-,t)=-\frac{t}{M^3}v^{(3)}(0)+\frac{1}{M^3}v^{(2)}(0)$
then \eqref{eq2.1:W3d2} rewrites
\begin{equation*}
\label{eq2.1:W3d2bis}
\left\{
\begin{aligned}
& W^{3/2}_t(w,t) - W_{ww}^{3/2}(w,t)=0, & (w,t)\in \mathbb{R}\times(0,T), \\
& \lim_{w\to +\infty} \biggl( W^{3/2}(w,t) -  \frac{y_0^{(3)}(0)}{6}(w^3+6wt) \biggr) =0, & t\in (0,T),\\
&  \lim_{w\to -\infty} \biggl( W^{3/2}(w,t) +  \frac{v^{(3)}(0)}{6M^3}(w^3+6wt) -  \frac{v^{(2)}(0)}{M^3} w\biggr)=0, & t\in (0,T).
\end{aligned}
\right.
\end{equation*}
A solution is given by $W^{3/2}(w,t)=(H(\cdot,t)\star g_{3/2})(w)$ with $g_{3/2}(s)= \frac{s^3}{6} y_0^{(3)}(0)$ if $s\geq 0$ and $g_{3/2}(s)=-\frac{s^3}{6} \frac{1}{M^3}v^{(3)}(0)+s\frac{1}{M^3}v^{(2)}(0) $ if $s<0$. 
Explicitly 
\begin{equation}
\label{def.W32}
\left\{
\begin{aligned}
&W^{3/2}(w,t)= \biggl(\frac{w^3}{2}+3tw\biggr)\biggl( erf\biggl(\frac{w}{2\sqrt{t}}\biggr)(h^+-h^-) + (h^++h^-)\biggr) \\
&\hspace{2cm} + (h^+-h^-)(4t+w^2)\sqrt{\frac{t}{\pi}} e^{-\frac{w2}{4t}}
-f^-\frac{\sqrt{t}}{\sqrt{\pi}} e^{-\frac{w^2}{4t}}  + \frac{f^-}{2}w \,erfc\biggl(\frac{w}{2\sqrt{t}}\biggr),\\
& f^-=\frac{v^{(2)}(0)}{M^3}, \quad h^{+}=\frac{y_0^{(3)}(0)}{6}, \quad h^{-}=-\frac{v^{(3)}(0)}{6M^3}.
\end{aligned}
\right.
\end{equation} 

\subsection{Composite asymptotic approximation outside of the boundary layer along $x=1$} 

One usual way to define a composite asymptotic approximation is to introduce a truncature $\mathcal{X}_{\eps}$ function in order to make the link between two expansions. This leads to an approximation of the form $\mathcal{X}_{\eps}(x,t)\sum_{k\geq 0} \eps^k y^{k} + (1-\mathcal{X}_{\eps}(x,t)) \sum_{k\geq 0} \eps^{k/2} W^{k/2}$. We refer to \cite{amirat_munch} where this strategy is employed. However this technique leads to tedious computations when the error analysis is performed. Instead, we define here a so-called composite approximation obtained by adding, at each order, the inner and outer expansions and then by subtracting their common part. 

At the first order, the common part of $y^0(x,t)$ and $W^0(w,t)$ (defined by \eqref{2.y0def} and \eqref{def.W0} respectively), is equal to $y_0(0)$ for $x>Mt$ and to $v(0)$ for $x<Mt$. Thus, the first term of the composite approximation outside the boundary layer along $x=1$ is given by
\begin{equation}
\label{def.p0}
p^0(x,t)=
\left\{
\begin{aligned}
& y_0(x-Mt)+ W^0(w,t)-y_0(0),  & x>Mt, \\
& v\biggl(t-\frac{x}{M}\biggr) + W^0(w,t)-v(0), & x\leq Mt. 
\end{aligned}
\right.
\end{equation}
Clearly, the function $p^0$ is continuous along the characteristic:
\begin{equation*}
lim_{x-Mt \to 0^{\pm}} p^0(x,t)=\frac{y_0(0)+v(0)}{2}.
\end{equation*} 

The second term of the composite approximation is given by
\begin{equation}
\label{def.p12}
p^{1/2}(x,t)=
\left\{
\begin{aligned}
& W^{1/2}(w,t)-y_0^{(1)}(0)\,w, &x> Mt, \\
& W^{1/2}(w,t)+\frac{1}{M}v^{(1)}(0)\, w, &x\leq Mt,
\end{aligned}
\right.
\end{equation}
where $W^{1/2}$ is defined by \eqref{def.W12}. Clearly, $p^{1/2}$
is also continuous along the characteristic. 

At the next order, we define the function
\begin{equation*}
p^1(x,t)=y^1(x,t)+W^{1}(w,t)-\biggl(\frac{w^2}{2}y_{xx}^0((Mt)^\pm,t)+y^1((Mt)^\pm,t)\biggr),
\end{equation*}
that is 
\begin{equation}
\label{def.p1}
p^1(x,t)=
\left\{
\begin{aligned}
& t y_0^{(2)}(x-Mt) + W^1(w,t)- \biggl(\frac{w^2+2t}{2}\biggr)y_0^{(2)}(0), &x\geq Mt, \\
& \frac{x}{M^3}v^{(2)}(t-\frac{x}{M})+W^1(w,t)-\biggl(\frac{w^2+2t}{2M^2}\biggr)v^{(2)}(0), &x\leq Mt, \\
\end{aligned}
\right.
\end{equation}
with $W^1$ given by \eqref{def.W1}. Eventually, we define $p^{3/2}$ as follows: 
\begin{equation*}
p^{3/2}(x,t)=W^{3/2}(w,t)-\biggl(\frac{w^3}{3 !} y^0_{xxx}((Mt)^\pm,t) +  w y_x^1( (Mt)^\pm,t)\biggr),
\end{equation*}
that is
\begin{equation}
\label{def.p32}
p^{3/2}(x,t)=
\left\{
\begin{aligned}
& W^{3/2}(w,t)-\biggl(\frac{w^3}{3 !} +  t w \biggr)y_0^{(3)}(0), &x\geq Mt, \\
& W^{3/2}(w,t)-\biggl(-\frac{w^3}{3 ! M^3} v^{(3)}(0) + \frac{w}{M^3}v^{(2)}(0)
-\frac{tw}{M^3}v^{(3)}(0)\biggr), &x\leq Mt, \\
\end{aligned}
\right.
\end{equation}
where $W^{3/2}$ is given by \eqref{def.W32}.

Then the following quantity is defined to be an asymptotic approximation of $y^{\eps}$, outside of the boundary layer along $x=1$, 
\begin{equation}
\label{def.pmeps}
p_3^{\eps}(x,t)=\sum_{k=0}^3 \eps^{k/2} p^{k/2}(x,t), \quad (x,t)\in Q_T. 
\end{equation}
We easily verify the following property. 

\begin{proposition}
\label{prop1}
Assume that $v\in C^3([0,T])$ and $y_0\in C^3([0,1])$. Then the functions 
 $(p^0+\sqrt{\eps}p^{1/2})$, $(p^1+\sqrt{\eps}p^{3/2})$ and $p_3^{\eps}$ belong to $C^1([0,1]\times]0,T])$.
\end{proposition}

\subsection{Inner expansion along $x=1$}\label{innerexpansionx1}

We now turn back to the construction of the inner expansion along $x=1$. 
Putting $\displaystyle \sum_{k=0}^m \varepsilon^{k/2} Y^{k/2}(z,\tau,t)$ into equation~$\eqref{eq:2.1}_1$, the identification of the powers of $\varepsilon$ yields
\begin{align*}
& \varepsilon^{0}: \quad Y^0_{zz}(z,\tau,t) + M Y_{z}^0(z,\tau,t)=0,  \\
& \varepsilon^{1/2}: \quad Y^{1/2}_{zz}(z,\tau,t) + M Y_{z}^{1/2}(z,\tau,t)=-Y_\tau^{0}(z,\tau,t) \\ 
& \varepsilon^{k/2}: \quad Y^{k/2}_{zz}(z,\tau,t) + M Y_{z}^{k/2}(z,\tau,t)=Y_t^{(k-2)/2}(z,\tau,t) -Y_\tau^{(k-1)/2}(z,\tau,t),  \quad \mbox{ for any } 2\leq k\leq m. 
\end{align*}

It is important to note that the functions $Y^k$ depends on three variables, namely $z, t$ but also $\tau$. As it is standart, the variable $z=(1-x)/\eps$ is introduced to describe the boundary layer at $x=1^-$. 
Here, the variable $\tau=(1/M-t)/\sqrt{\eps}$ allows to take into account the interaction of the internal and boundary layer.  If we do not introduce this variable $\tau$, we see notably that $Y^{1/2}$ solves the same ordinary differential equation than $Y^0$, and the analysis (detailed in Section \ref{section_convergence}) leads to an error estimate (for the $L^\infty(L^2)$ norm) of the order $\sqrt{\eps}$ only, independently of the number of terms in the various expansions. This point, also described in \cite{howes1989} to discuss the problem (\ref{Pb_elliptic}), is therefore crucial. 

We impose $Y^{k/2}(0,\tau)=0$ for any $0 \leq k\leq m$. To get the asymptotic matching conditions we write that, for any fixed $\tau$ and large $z$,
\begin{align*}
& Y^0(z,\tau,t) + \varepsilon^{1/2} Y^{1/2}(z,\tau,t) + \varepsilon Y^1(z,\tau,t) + \cdots + \varepsilon^{m/2} Y^{m/2}(z,\tau,t) \\
& =p^0(x,t) + \varepsilon^{1/2} p^{1/2}(x,t) + \varepsilon^1 p^1(x,t) + \cdots + \varepsilon^{m/2} p^{m/2}(x,t) + \mathcal{O}(\varepsilon^{m+1}).
\end{align*}

In order to identify at each order the appropriate matching conditions, we need to rewrite the right-hand side of the above equality in terms of $z$ and $\tau$, $t$ being fixed. 
Writing that $x=1-\eps z$, $w=M\tau-\sqrt{\eps}z$, and using Taylor expansions, we have successively
$$
\begin{aligned}
p^0(x,t)& =y^0(x,t)+W^0(w,t)-y^0((Mt)^\pm,t)= y^0(1-\eps z,t)+ W^0(M\tau-\sqrt{\eps}z,t)-y^0((Mt)^\pm,t)\\
&=y^0(1,t)+W^0(M\tau,t)-y^0((Mt)^\pm,t) -\sqrt{\eps}z W^0_w(M\tau,t)+\eps \biggl(-z (y^0)_x(1,t)+ \frac{z^2}{2}W_{ww}^0(M\tau,t)\biggr)\\
& \quad -\eps^{3/2} \frac{z^3}{6}W^0_{www}(M\tau,t)+\mathcal{O}(\eps^{2}),
\end{aligned}
$$
$$ 
\begin{aligned}
\sqrt{\eps}p^{1/2}(x,t) =&\eps^{1/2}\biggl(W^{1/2}(w,t)-w (y^0)_x((Mt)^\pm,t)\biggr)\\
=&\eps^{1/2}\biggl(W^{1/2}(M\tau-\sqrt{\eps}z,t)-(M\tau-\sqrt{\eps}z) (y^0)_x((Mt)^\pm,t) \biggr)  \\
=&\eps^{1/2}\biggl(W^{1/2}(M\tau,t)-M\tau (y^0)_x((Mt)^\pm,t)\biggr)-\eps z \biggl(W_w^{1/2}(M\tau,t)-(y^0)_x((Mt)^\pm,t)\biggr)\\
&+\eps^{3/2}\frac{z^2}{2}W_{ww}^{1/2}(M\tau,t)+\mathcal{O}(\eps^{2}),
\end{aligned}
$$
$$ 
\begin{aligned}
\eps p^1(x,t)=& \eps \biggl(y^1(x,t)+ W^1(w,t)-\biggl(\frac{w^2}{2}y_{xx}^0((Mt)^\pm,t)+y^1((Mt)^\pm,t)\biggr)\biggr)\\
=&\eps \biggl(y^1(1-\eps z,t)+ W^1(M\tau-\sqrt{\eps}z,t)-\biggl(\frac{(M\tau-\sqrt{\eps}z)^2}{2}y_{xx}^0((Mt)^\pm,t)+y^1((Mt)^\pm,t)\biggr)\biggr) \\
=& \eps\biggl(y^1(1,t)+W^1(M\tau,t)-\biggl(\frac{M^2\tau^2}{2}y_{xx}^0((Mt)^\pm,t)+y^1((Mt)^\pm,t)\biggr)\biggr)\\
&+\eps^{3/2}z\biggl(-W^{1}_{w}(M\tau,t)+M\tau  y_{xx}^0((Mt)^\pm,t)\biggr)+\mathcal{O}(\eps^{2}),
\end{aligned}
$$
and
$$ 
\begin{aligned}
& \eps^{3/2} p^{3/2}(x,t)= \eps^{3/2}\left(
W^{3/2}(w,t)-\biggl(\frac{w^3}{3 !} y^0_{xxx}((Mt)^\pm,t) +  w y_x^1( (Mt)^\pm,t)\biggr)
\right) \\
& \quad\quad = \eps^{3/2}\left(W^{3/2}(M\tau-\sqrt{\eps}z,t)-\left(
\frac{(M\tau-\sqrt{\eps}z)^3}{3 !} y^0_{xxx}((Mt)^\pm,t)
+ (M\tau-\sqrt{\eps}z) y_x^1( (Mt)^\pm,t)
\right)\right)\\
& \quad\quad = \eps^{3/2}\left(W^{3/2}(M\tau,t)-\left(
\frac{M^3\tau^3}{3 !} y^0_{xxx}((Mt)^\pm,t)
+ M\tau y_x^1( (Mt)^\pm,t)
\right)\right) +\mathcal{O}(\eps^{2}).
\end{aligned}
$$
We deduce that 
$$
\sum_{k=0}^3 \eps^{k/2} p^{k/2}(x,t) = C_0(z,\tau,t)+\eps^{1/2} C_{1/2}(z,\tau,t)+\eps C_1(z,\tau,t)+\eps^{3/2} C_{3/2}(z,\tau,t)+\mathcal{O}(\eps^2),
$$
with 
\begin{equation}
\label{Cztaut}
\begin{aligned}
 C_0(z,\tau,t)=&y^0(1,t)+W^0(M\tau,t)-y^0((Mt)^\pm,t), \\
 C_{1/2}(z,\tau,t)=&W^{1/2}(M\tau,t)-M\tau (y^0)_x((Mt)^\pm,t)-z W^0_w(M\tau,t),\\
 C_1(z,\tau,t)=&y^1(1,t)+W^1(M\tau,t)-\biggl(\frac{M^2\tau^2}{2}y_{xx}^0((Mt)^\pm,t)+y^1((Mt)^\pm,t)\biggr)\\
 & +z \biggl(-W_w^{1/2}(M\tau,t)+(y^0)_x((Mt)^\pm,t)- (y^0)_x(1,t)\biggr)+ \frac{z^2}{2}W_{ww}^0(M\tau,t),\\
 C_{3/2}(z,\tau,t)=& p^{3/2}(1,t)+ z\biggl(-W^{1}_{w}(M\tau,t)+M\tau  y_{xx}^0((Mt)^\pm,t)\biggr)+\frac{z^2}{2}W_{ww}^{1/2}(M\tau,t)-\frac{z^3}{6}W^0_{www}(M\tau,t).
\end{aligned}
\end{equation}
Therefore, the matching conditions read 
\begin{equation}
Y^{k/2}(z,\tau,t) \sim C_{k/2}(z,\tau,t), \quad \mbox{ as } z\to \infty,  \quad \mbox{ for } k=0,\dots,3. \label{matching_condition_z}
\end{equation}
In the sequel, we use the notations
$$
C_{k/2}(z,\tau,t)=\sum_{p=0}^{k} C_{k/2,p}(\tau,t) z^p.
$$

We define $Y^0$ as a solution of
\begin{equation*}
\label{eq2.1:Y0}
\left\{
\begin{aligned}
& Y^0_{zz}(z,\tau,t) + M Y_{z}^0(z,\tau,t)=0, & (z,\tau,t)\in \mathbb{R}^\star_+\times \mathbb{R}\times (0,T), \\
& Y^0(0,\tau,t)=0, \quad  \lim_{z\to +\infty} Y^0(z,\tau,t)=C_{0,0}(\tau,t), & t\in (0,T). 
\end{aligned}
\right.
\end{equation*}
The solution is 
\begin{equation}
\label{2.Y0.def}
Y^0(z,\tau,t)=C_{0,0}(\tau,t)\left(1-e^{-Mz}\right), \quad (z,\tau,t)\in \mathbb{R}_+\times\mathbb{R}\times [0,T],
\end{equation}
with 
\begin{equation*}
\label{C00.def}
C_{0,0}(\tau,t)=
\left\{
\begin{aligned}
& y_0(1-Mt)+W^0(M\tau,t) - y_0(0), & 1\geq Mt,\\
& v\biggl(t-\frac{1}{M}\biggr)+W^0(M\tau,t)-v(0), & 1 \leq Mt.
\end{aligned}
\right.
\end{equation*}

Similarly, we define $Y^{1/2}$ as a solution of 
\begin{equation*}
\label{eq2.1:Y1s2}
\left\{
\begin{aligned}
& Y^{1/2}_{zz}(z,\tau,t) + M Y_{z}^{1/2}(z,\tau,t)=-Y_\tau^0(z,\tau,t), & (z,\tau,t)\in \mathbb{R}^\star_+\times \mathbb{R}\times (0,T), \\
& Y^{1/2}(0,\tau,t)=0, \quad  \lim_{z\to +\infty} (Y^{1/2}(z,\tau,t)-C_{1/2}(z,\tau,t))=0, & t\in (0,T). 
\end{aligned}
\right.
\end{equation*}
Writing that $Y_\tau^0(z,\tau,t)=C_{0,0,\tau}(1-e^{-Mz})$, we obtain, for 
$(z,\tau,t)\in \mathbb{R}^\star\times \mathbb{R}\times (0,T)
$,~\footnote{
The general solution of the ordinary differential equation $Y_{zz}+MY_z= (\alpha+\beta z + \frac{\gamma}{2} z^2) + e^{-Mz} (-\alpha+\beta z - \frac{\gamma}{2} z^2), \,Y(0)=0$ is in the form 
$$
Y(z)= \biggl(A+Bz+C\frac{z^2}{2}+D\frac{z^3}{6}\biggr)+e^{-Mz}\biggl(-A+Bz-C\frac{z^2}{2}+D\frac{z^3}{6}\biggr),
$$
with 
$$
A =\frac{1}{M}(K-B), \quad B=  \frac{1}{M}(\alpha-C), \quad C=\frac{1}{M}(\beta-D), \quad D= \frac{\gamma}{M}, \quad K
\mbox{ arbitrary constant. }
$$
}
\begin{equation}
\label{2.Y12.def}
\begin{aligned}
Y^{1/2}(z,\tau,t)=& \biggl(C_{1/2,0}(\tau,t)-\frac{C_{0,0,\tau}(\tau,t)}{M}z\biggr)+ e^{-Mz}\biggl(-C_{1/2,0}(\tau,t)-\frac{C_{0,0,\tau}(\tau,t)}{M}z\biggr), \\
=&\biggl(C_{1/2,0}(\tau,t)-W^0_w(M\tau,t)z\biggr)+ e^{-Mz}\biggl(-C_{1/2,0}(\tau,t)-W^0_w(M\tau,t)z\biggr).
\end{aligned}
\end{equation}

Then, we define $Y^1$ as a solution of 
\begin{equation*}
\label{eq2.1:Y1}
\left\{
\begin{aligned}
& Y^1_{zz}(z,\tau,t) + M Y_{z}^1(z,\tau,t)=Y^0_t(z,\tau,t)-Y_\tau^{1/2}(z,\tau,t), & (z,\tau,t)\in \mathbb{R}^\star\times \mathbb{R}\times (0,T), \\
& Y^1(0,\tau,t)=0, \quad  \lim_{z\to +\infty} (Y^1(z,\tau,t)-C_1(z,\tau,t))=0, & t\in (0,T). 
\end{aligned}
\right.
\end{equation*}
We have
$$Y_\tau^{1/2}(z,\tau,t)=\biggl(C_{1/2,0,\tau}(\tau,t)-MW_{ww}^0(M\tau,t)z\biggr)+ e^{-Mz}\biggl(-C_{1/2,0,\tau}(\tau,t)-M
W_{ww}^0(M\tau,t)z\biggr),
$$ 
and 
$$
Y_t^0(z,\tau,t)=C_{0,0,t}(\tau,t)(1-e^{-Mz})=\biggl(y^0_t(1,t)+W^0_t(M\tau,t)\biggr)(1-e^{-Mz}),
$$
then 
$$
\begin{aligned}
Y^0_t(z,\tau,t)-Y_\tau^{1/2}(z,\tau,t)=&\biggl(C_{0,0,t}(\tau,t)-C_{1/2,\tau}(0,\tau,t)+MW_{ww}^0(M\tau,t)z\biggr)\\
&+ e^{-Mz}\biggl(-C_{0,0,t}(\tau,t)+C_{1/2,\tau}(0,\tau,t)+MW_{ww}^0(M\tau,t)z\biggr).
\end{aligned}
$$
We have
\begin{equation*}
C_{0,0,t}(\tau,t)-C_{1/2,\tau}(0,\tau,t)=
\left\{
\begin{aligned}
& -My_0^{(1)}(1-Mt) + W^0_t(M\tau,t)
-M W^{1/2}_w(M\tau,t) + M y_0^{(1)}(0), & 1> Mt,\\
& v^{(1)}(t-\frac{1}{M}) + W^0_t(M\tau,t) -M W^{1/2}_w(M\tau,t) -v^{(1)}(0), & 1< Mt,
\end{aligned}
\right.
\end{equation*}
and note that $C_{0,0,t}(\tau,t)-C_{1/2,\tau}(0,\tau,t)$ has no jump.
Explicitly, we obtain
\begin{equation}
\label{2.Y1.def}
Y^1(z,\tau,t)=\left(A+Bz+C\frac{z^2}{2}\right)+e^{-Mz}\left(-A+Bz-C\frac{z^2}{2}\right),
\end{equation}
with
\begin{equation}
\label{ABC.Y1.def}
\left\{
\begin{aligned}
A=& C_1(0,\tau,t)=y^1(1,t)+W^1(M\tau,t)-\biggl(\frac{M^2\tau^2}{2}y_{xx}^0((Mt)^\pm,t)+y^1((Mt)^\pm,t)\biggr),\\
B=&\frac{1}{M}\biggl(C_{0,0,t}(\tau,t)-C_{1/2,0,\tau}(\tau,t)-W_{ww}^0(M\tau,t)\biggr)=-y_x^0(1,t)-W_w^{1/2}(M\tau,t)+y_x^0((Mt)^\pm,t), \\
C=& W_{ww}^0(M\tau,t).
\end{aligned}
\right.
\end{equation}

Similarly, we define $Y^{3/2}$ as a solution of 
\begin{equation*}
\label{eq2.1:Y32}
\left\{
\begin{aligned}
& Y^{3/2}_{zz}(z,\tau,t) + M Y_{z}^{3/2}(z,\tau,t)=Y^{1/2}_t(z,\tau,t)-Y_\tau^1(z,\tau,t), & (z,\tau,t)\in \mathbb{R}^\star_+\times\mathbb{R}^\star_+\times (0,T), \\
& Y^{3/2}(0,\tau,t)=0, \quad  \lim_{z\to +\infty} (Y^{3/2}(z,\tau,t)-C_{3/2}(z,\tau,t))=0, & t\in (0,T). 
\end{aligned}
\right.
\end{equation*}
We have 
\begin{equation*}
Y^{1/2}_t(z,\tau,t)=\biggl(C_{1/2,t}(0,\tau,t)-W_{wt}^0(M\tau,t)z\biggr)+ e^{-Mz}\biggl(-C_{1/2,t}(0,\tau,t)-W_{wt}^0(M\tau,t)z\biggr),
\end{equation*}
and 
\begin{equation*}
\label{Y^1}
\begin{aligned}
Y^1_\tau(z,\tau,t)=&M\biggl(W^1_w(M\tau,t)-M \tau y^0_{xx}((Mt)^\pm,t)-W^{1/2}_{ww}(M\tau,t)z+W^0_{www}(M\tau,t)\frac{z^2}{2}\biggr)\\
&+Me^{-Mz}\biggl(-W^1_w(M,\tau,t)+M \tau y^0_{xx}((Mt)^\pm,t)-W^{1/2}_{ww}(M\tau,t)z- W^0_{www}(M\tau,t)\frac{z^2}{2}\biggr).
\end{aligned}
\end{equation*}
Explicitly
\begin{equation}
\label{2.Y32.def}
Y^{3/2}(z,\tau,t)=\left(\widetilde{A}+\widetilde{B}z+\widetilde{C}\frac{z^2}{2}+\widetilde{D}\frac{z^3}{6}\right)+e^{-Mz}\left(-\widetilde{A}+\widetilde{B}z-\widetilde{C}\frac{z^2}{2}+\widetilde{D}\frac{z^3}{6}\right),
\end{equation}
with
\begin{equation}
\label{ABC.Y32.def}
\left\{
\begin{aligned}
\widetilde{A}=& C_{3/2}(0,\tau,t)=W^{3/2}(M\tau,t)-\biggl(\frac{(M\tau)^3}{3 !} y^0_{xxx}((Mt)^\pm,t) +  (M\tau) y_x^1( (Mt)^\pm,t)\biggr), \\
\widetilde{B}=&\frac{1}{M}\biggl(C_{1/2,t}(0,\tau,t)-M\biggl(W_w^1(M\tau,t)-M\tau y_{xx}^0((Mt)^\pm,t)\biggr)-\widetilde{C}\biggr) \\
\; = & -\biggl(W_w^1(M\tau,t)-M\tau y_{xx}^0((Mt)^\pm,t)\biggr), \\
\widetilde{C}=& W_{ww}^{1/2}(M\tau,t), \quad \widetilde{D}= -W_{www}^0(M\tau,t).\\
\end{aligned}
\right.
\end{equation}

\subsection{Asymptotic composite approximation in $Q_T$}

We are now in position to define what is supposed to be an asymptotic approximation of the solution $y^{\eps}$. 
We proceed as before by adding at each order the term $p^{k/2}$, approximation outside of the boundary layer along $x=1$, and the term $Y^{k/2}$, approximation in the boundary layer along $x=1$, then subtracting their common part.  A the first order, the composite approximation is given by
\begin{equation}
\label{def.P0}
\begin{aligned}
P^0(x,t) & = p^0(x,t)+Y^0(z,\tau,t) -C_0(\tau,t)=p^0(x,t)-C_0(\tau,t)e^{-Mz}\\
& =y^0(x,t)+W^0(w,t)-y^0((Mt)^\pm,t)
-\left(y^0(1,t)+W^0(M\tau,t)-y^0((Mt)^\pm,t)\right) e^{-Mz}.
\end{aligned}
\end{equation}
so that $P^0(1,t)=0$ for all $t\geq 0$. Repeating the process, we define 
\begin{equation}
\label{def.Pk2}
P^{k/2}(x,t)=p^{k/2}(x,t)+Y^{k/2}(z,\tau,t) -C_{k/2}(z,\tau,t), \quad 1\leq k \leq 3,
\end{equation}
then we define an asymptotic composite approximation of $y^\eps$ in $Q_T$ by
\begin{equation}
\label{def.P3eps}
P^{\eps}(x,t)=\sum_{k=0}^3 \eps^{k/2} P^{k/2}(x,t), \quad (x,t)\in Q_T.
\end{equation}

\section{The sequence $(\widetilde{P}^{\eps})$: another sequence of approximate solutions}\label{section3}

As we will see in the next section, approximation $P^{\eps}$ leads to an error estimate of order $\sqrt{\eps}$ for the $L^\infty(L^2)$ norm, which is not satisfactory. In this section, we construct another sequence of approximate solutions. The construction is made by modifying the previous sequence $(P^{\eps})$ defined by \eqref{def.P3eps}. Precisely, we replace the function $W^0$, see (\ref{def.W0}), introduced in Section \ref{innerexp}, to describe the inner layer along the characteristic $\{(x,t)\in Q_T,x-Mt=0\}$. The reason of this modification is the following one. Recalling that  $w=(x-Mt)/\sqrt{\eps}$, 
for all $\eps>0$ and $x>0$, we check that $W^0\big(\frac{x}{\sqrt{\eps}},0\big)=y_0(0)$. On the other hand, $\lim_{t\to 0^+} W^0\big(\frac{-Mt}{\sqrt{\eps}},t\big)=\frac{y_0(0)+v(0)}{2}$, different from $v(0)$ when $v(0)\neq y^0(0)$.
For this reason, the use of this $W^0$ function generates an artificial boundary layer in the approximate solution along the line $x=0$, above the characteristic in the neighborhood of $t=0$. This boundary layer propagates inside the domain, in the neighborhood of the characteristic and affect the quality of the approximation. To avoid this fact, we consider, instead of the function $g_0$ defined in (\ref{def.g0}), the following second choice: 
\begin{equation*}
g^{\eps}_0(w)=\left\{
\begin{aligned}
& y_0(0), & w\geq 0,\\
& v(0)+(v(0)-y_0(0))e^\frac{Mw}{\sqrt{\eps}}, & w<0,
\end{aligned}
\right.
\end{equation*} 
leading to (recalling that $H(w,t)=(4\pi t)^{-1/2}e^{-w^2/(2t)}$)
\begin{equation}
\label{def.W0eps}
W^0_\eps(w,t)=(H(\cdot,t)\star g_0^\eps)(w)=W^0(w,t)+ U^0_\eps(w,t),
\end{equation}
where
\begin{equation}
\label{def.U0eps}
U^0_\eps(w,t)=\frac{v(0)-y_0(0)}{2}e^{\frac{Mw}{\sqrt{\eps}}+\frac{M^2t}{\eps}}erfc\left(\frac{w}{2\sqrt{t}}+\frac{M\sqrt{t}}{\sqrt{\eps}}\right).
\end{equation}
In particular, we have remarkably
\begin{equation}
W_\eps^0\biggl(-\frac{Mt}{\sqrt{\eps}},t\biggr)-v(0)=0, \quad \forall t\in (0,T] \label{propertyW0eps}
\end{equation}
and still the property $W_\eps^0(\frac{x}{\sqrt{\eps}},0)=y_0(0)$ for all $x>0$. Actually, the function $\widetilde{W}^0_{\eps}(x,t)=W^0_{\eps}(\frac{x-Mt}{\sqrt{\eps}},t)$ solves the equation (we refer notably to \cite{Shih07})
\begin{equation*}
\label{eq:transportinfini}
	\left\{
		\begin{aligned}
   			& \widetilde{W}^0_{\eps,t} + M \widetilde{W}^0_{\eps,x} - \eps \widetilde{W}^0_{\eps,xx}=0, 				& (x,t)\in \mathbb{R}^+\times \mathbb{R}^+,    \\
   			& \widetilde{W}^0_{\eps}(0,t)  = v(0), 	  		& t\in  \mathbb{R}^+, \\
   			& \widetilde{W}^0_{\eps}(x, 0) = y_0(0), 						  		& x\in  \mathbb{R}^+.
   		\end{aligned} 
 	\right.
\end{equation*}
$\widetilde{W}^0_{\eps}$ is therefore the explicit solution of the advection-diffusion equation defined over $(x,t)\in \mathbb{R}^\star\times \mathbb{R}^{\star}$ associated to constant initial and boundary data.
Remark that on the line $x=0$, the function $U_\eps^0$ satisfies 
$$
U_\eps^0\biggl(-\frac{Mt}{\sqrt{\eps}},t\biggr)= \frac{v(0)-y_0(0)}{2}  erfc\biggl(\frac{M\sqrt{t}}{2\sqrt{\eps}}\biggr)\leq  \frac{v(0)-y_0(0)}{2}  e^{-\frac{M^2t}{4\eps}}, \quad\forall t\geq 0
$$
and gets concentrated as $\eps$ goes to zero in the neighborhood of $t=0$ and allows to eliminate the artificial boundary layer mentioned above.

We now choose $W_\eps^0(w,t)$ (defined by \eqref{def.W0eps}, \eqref{def.U0eps}), instead of $W^0(w,t)$, as a solution of equation \eqref{eq2.1:W0}. By analogy with \eqref{def.p0} we define the function
\begin{equation}
\label{def.p0tilde}
p_\eps^0(x,t)=
\left\{
\begin{aligned}
& y_0(x-Mt)+ W_\eps^0(w,t)-y_0(0),  & x>Mt, \\
& v\biggl(t-\frac{x}{M}\biggr) + W_\eps^0(w,t)-v(0), & x\leq Mt. 
\end{aligned}
\right.
\end{equation}
We note that $p_\eps^0$ is continuous along the characteristic:
\begin{equation*}
lim_{x-Mt \to 0^{\pm}} p_\eps^0(x,t)=\frac{y_0(0)+v(0)}{2} + \frac{v(0)-y_0(0)}{2}e^{\frac{M^2t}{\eps}}erfc\left(\frac{M\sqrt{t}}{\sqrt{\eps}}\right).
\end{equation*} 
The following quantity is defined to be an asymptotic approximation of $y^{\eps}$, outside the boundary layer along $x=1$,
\begin{equation}
\label{def.pmepstilde}
\widetilde{p}_3^{\eps}(x,t)=p_\eps^0(x,t) + \sum_{k=1}^3 \eps^{k/2} p^{k/2}(x,t), \quad (x,t)\in Q_T, 
\end{equation}
where the functions $p^{k/2}$, $k=1,2,3$, are defined by \eqref{def.p12}, \eqref{def.p1}, \eqref{def.p32}, respectively.
We have the analog of Proposition \ref{prop1}.

\begin{proposition}
\label{prop2}
Assume that $v\in C^3([0,T])$ and $y_0\in C^3([0,1])$. Then the functions 
$(p^0_\eps+\sqrt{\eps}p^{1/2})$, and $\widetilde{p}_3^{\eps}$ belong to $C^1([0,1]\times]0,T])$.
\end{proposition}

Let us now introduce the notations:
\begin{equation}
\label{Cztauteps}
\begin{aligned}
 C^\eps_0(z,\tau,t)=&y^0(1,t)+W_\eps^0(M\tau,t)-y^0((Mt)^\pm,t), \\
 C^\eps_{1/2}(z,\tau,t)=&W^{1/2}(M\tau,t)-M\tau (y^0)_x((Mt)^\pm,t)-z W^0_{\eps,w}(M\tau,t),\\
 C^\eps_1(z,\tau,t)=&y^1(1,t)+W^1(M\tau,t)-\biggl(\frac{M^2\tau^2}{2}y_{xx}^0((Mt)^\pm,t)+y^1((Mt)^\pm,t)\biggr)\\
 & +z \biggl(-W_w^{1/2}(M\tau,t)+(y^0)_x((Mt)^\pm,t)- (y^0)_x(1,t)\biggr)+ \frac{z^2}{2}W_{\eps,ww}^0(M\tau,t),\\
 C^\eps_{3/2}(z,\tau,t)=& p^{3/2}(1,t)+ z\biggl(-W^{1}_{w}(M\tau,t)+M\tau  y_{xx}^0((Mt)^\pm,t)\biggr)+\frac{z^2}{2}W_{ww}^{1/2}(M\tau,t)-\frac{z^3}{6}W^0_{\eps,www}(M\tau,t).
\end{aligned}
\end{equation}
Continuing the construction, we define the analog of the function $Y^0$ by
\begin{equation}
\label{2.Y0eps.def}
Y_\eps^0(z,\tau,t)=C^\eps_{0}(z,\tau,t)\left(1-e^{-Mz}\right), \quad (z,\tau,t)\in \mathbb{R}_+\times\mathbb{R}\times [0,T],
\end{equation}
%
%
The analog of $Y^{1/2}$ is the function
\begin{equation}
\label{2.Y12eps.def}
Y_\eps^{1/2}(z,\tau,t)
=\biggl(C_{1/2,0}(\tau,t)-W^0_{\eps,w}(M\tau,t)z\biggr)+ e^{-Mz}\biggl(-C_{1/2,0}(\tau,t)-W^0_{\eps,w}(M\tau,t)z\biggr),
\end{equation}
where $C_{1/2,0}(\tau,t)$ is defined by \eqref{Cztaut}. The analog of $Y^1$ is given by
\begin{equation}
\label{2.Y1eps.def}
Y_\eps^1(z,\tau,t)=\left(A+Bz+C^\eps\frac{z^2}{2}\right)+e^{-Mz}\left(-A+Bz-C^\eps\frac{z^2}{2}\right),
\end{equation}
with $A$, $B$ defined by $(\ref{ABC.Y1.def})_1$ and $(\ref{ABC.Y1.def})_2$, respectively, and 
$C^\eps=W_{\eps,ww}^0(M\tau,t)$. Eventually, the analog of $Y^{3/2}$ is given by
\begin{equation}
\label{2.Y32eps.def}
Y_\eps^{3/2}(z,\tau,t) = \left(\widetilde{A}+\widetilde{B}z+\widetilde{C}\frac{z^2}{2}+\widetilde{D}^\eps\frac{z^3}{6}\right)+e^{-Mz}\left(-\widetilde{A}+\widetilde{B}z-\widetilde{C}\frac{z^2}{2}+\widetilde{D}^\eps\frac{z^3}{6}\right),
\end{equation}
with $\widetilde{A}$, $\widetilde{B}$, $\widetilde{C}$, defined by $(\ref{ABC.Y32.def})_1$, $(\ref{ABC.Y32.def})_2$, $(\ref{ABC.Y32.def})_3$, and $\widetilde{D}^\eps= -W_{\eps,www}^0(M\tau,t).$

\begin{remark}
\label{remarkYkeps}
We check that the functions $Y^{k/2}_\eps$, $k=0,1,2,3,$ solve, respectively, the following equations:
 \begin{equation*}
\label{eq2.1:Y0eps}
\left\{
\begin{aligned}
& Y^0_{\eps,zz}(z,\tau,t) + M Y_{\eps,z}^0(z,\tau,t)=0, & (z,\tau,t)\in \mathbb{R}^\star_+\times \mathbb{R}\times (0,T), \\
& Y^0_{\eps}(0,\tau,t)=0, \quad  \lim_{z\to +\infty} Y^0_{\eps}(z,\tau,t)=C_0^\eps(0,\tau,t), & t\in (0,T); 
\end{aligned}
\right.
\end{equation*}

\begin{equation*}
\label{eq2.1:Y1s2eps}
\left\{
\begin{aligned}
& Y^{1/2}_{\eps,zz}(z,\tau,t) + M Y_{\eps,z}^{1/2}(z,\tau,t)=-Y_{\eps,\tau}^0(z,\tau,t), & (z,\tau,t)\in \mathbb{R}^\star_+\times \mathbb{R}\times (0,T), \\
& Y_{\eps}^{1/2}(0,\tau,t)=0, \quad  \lim_{z\to +\infty} (Y_{\eps}^{1/2}(z,\tau,t)-C_{1/2}^\eps(z,\tau,t))=0, & t\in (0,T);
\end{aligned}
\right.
\end{equation*}

\begin{equation*}
\label{eq2.1:Y1eps}
\left\{
\begin{aligned}
& Y^1_{\eps,zz}(z,\tau,t) + M Y_{\eps,z}^1(z,\tau,t)=Y^0_{\eps,t}(z,\tau,t)-Y_{\eps,\tau}^{1/2}(z,\tau,t), & (z,\tau,t)\in \mathbb{R}^\star\times \mathbb{R}\times (0,T), \\
& Y_\eps^1(0,\tau,t)=0, \quad  \lim_{z\to +\infty} (Y_\eps^1(z,\tau,t)-C_1^\eps(z,\tau,t))=0, & t\in (0,T); 
\end{aligned}
\right.
\end{equation*}

\begin{equation*}
\label{eq2.1:Y32eps}
\left\{
\begin{aligned}
& Y^{3/2}_{\eps,zz}(z,\tau,t) + M Y_{\eps,z}^{3/2}(z,\tau,t)=Y^{1/2}_{\eps,t}(z,\tau,t)-Y_{\eps,\tau}^1(z,\tau,t), & (z,\tau,t)\in \mathbb{R}^\star_+\times\mathbb{R}^\star_+\times (0,T), \\
& Y_\eps^{3/2}(0,\tau,t)=0, \quad  \lim_{z\to +\infty} (Y_\eps^{3/2}(z,\tau,t)-C_{3/2}^\eps(z,\tau,t))=0, & t\in (0,T). 
\end{aligned}
\right.
\end{equation*}
\end{remark}

We now define the first term of the asymptotic approximation in $Q_T$ by
\begin{equation}
\label{def.P0eps}
\begin{aligned}
P_\eps^0(x,t) & = p_\eps^0(x,t)+ Y_\eps^0(z,\tau,t)-C^\eps_0(z,\tau,t)e^{-Mz}\\
& = y^0(x,t)+W_\eps^0(w,t)-y^0((Mt)^\pm,t)
-\left(y^0(1,t)+W_\eps^0(M\tau,t)-y^0((Mt)^\pm,t)\right) e^{-Mz}.
\end{aligned}
\end{equation}
Then, setting 
\begin{equation}
\label{def.Pk2eps}
P_\eps^{k/2}(x,t)=p^{k/2}(x,t)+Y_\eps^{k/2}(z,\tau,t) -C^\eps_{k/2}(z,\tau,t), \quad 1\leq k \leq 3,
\end{equation}
the function $\widetilde{P}^{\eps}$ defined to be an asymptotic approximation of $y^\eps$ in $Q_T$ is given by
\begin{equation}
\label{def.P3epstilde}
\widetilde{P}^{\eps}(x,t)=\sum_{k=0}^3 \eps^{k/2} P_\eps^{k/2}(x,t), \quad (x,t)\in Q_T.
\end{equation}

\begin{remark}
\label{explicit}
It should be noted that, as for the sequence $(P^\eps)$ constructed in Section \ref{section2}, the sequence $(\widetilde{P}^{\eps})$ is constructed by using explicit formulae. These asymptotic approximations are very useful, for instance, from a numerical viewpoint (easily computable because using explicit formulae). Comparing with the sequence $({P}^{\eps})$ we will see that $(\widetilde{P}^{\eps})$ is a more accurate approximation of the solution $y^\eps$ of \eqref{eq:2.1}.
\end{remark}

In the next section we investigate the convergence of the sequence $(\widetilde{P}^{\eps})$. 

\section{Convergence of the sequence $(\widetilde{P}^{\eps})_{(\eps>0)}$. Rate of convergence}\label{section_convergence}

This section is devoted to a study on the convergence of the sequence $(\widetilde{P}^{\eps})$. 
We make the following assumptions:
\begin{equation}
\label{assumption}
y_0 \in C^4([0,1]), \quad v\in C^4([0,T]).
\end{equation}
Our main goal is to establish the following result.

\begin{theorem}
\label{convergence}
Let $y^{\eps}$ be the solution of \eqref{eq:2.1} and $\widetilde{P}^{\eps}$ the function defined by \eqref{def.P3epstilde}. Assume \eqref{assumption}. Then there exists two constants $c$ and $\eps_0$, $c$ independent of $\eps$, such that, for $0<\eps< \eps_0$, and any $\gamma\in (0,1/2]$
\begin{equation}
\label{rate}
 \Vert \widetilde{P}^{\eps}(\cdot,t)-y^{\eps}(\cdot,t) \Vert_{L^2(0,1)} \leq c\,\eps^{3/2} + c\,\eps^{1/2}
e^{-\frac{M^2}{2\eps^\gamma}t} \quad \forall t\in [0,T].
\end{equation}
\end{theorem}


In order to prove this theorem we need to establish a number of preliminary results. We define the error as follows: 
\begin{equation}
z^{\eps}(x,t)=\widetilde{P}^{\eps}(x,t)-y^{\eps}(x,t)-\theta^\eps(x,t), \quad (x,t)\in Q_T, \label{defzeps}
\end{equation}
where $\theta^\varepsilon$ is the initial layer corrector defined as a solution of the equation
\begin{equation}
\label{inlayer1}
\left\{
\begin{aligned}
& \theta^{\varepsilon}_t + M \theta^{\varepsilon}_x - \varepsilon \theta^{\varepsilon}_{xx} =0, & (x,t) \in Q_T,\\
& \theta^{\varepsilon}(0,t)=\theta^{\varepsilon}(1,t) =0, & t\in (0,T),\\
& \theta^{\varepsilon}(x,0)=\tilde{P}^{\eps}(x,0)-y^\eps(x,0), & x\in (0,1).
\end{aligned}
\right.
\end{equation}

The occurence of the initial layer is due to the fact that, by construction, the approximation $\widetilde{P}^{\eps}(1,t)$ vanishes for all $t$, including $t=0$, while the value of $y_0$ at $x=1$ may be different from zero. This introduces an error of order one but which get concentrated as $\eps$ decreases in the neighborhood of the point $(1,0)$.  Thank to the transport term, which push the solution to the right, it turns out that this error is damped out exponentially fast at $t$ increases. Introducing the auxiliary variable $\tau_1=t/\sqrt{\eps}$ coupled with the variable $z=(1-x)/\eps$, we may approximate precisely this so-called corner layer. The computations are very similar to the one performed to treat the interaction of the internal and boundary layers, in the neighborhood of $(1,1/M)$, by introducing the zoom variables $z$ and $\tau$ (see Section \ref{innerexpansionx1}). We refer to \cite[Section 4.1]{ColeBook81} where this is discussed. We emphasize however, that in view of the distribution of $\theta^{\eps}(\cdot,0)$ along $(0,1)$, this is not necessary for our objective to get an $L^\infty(L^2)$ estimate.

\subsection{Preliminary results}

\subsubsection{Estimate of the initial layer corrector $\theta^\eps$}

The following lemma gives an exponential decay property of the initial layer corrector.

\begin{lemma}
\label{prelim1}
Let $\theta^{\varepsilon}$ be the solution of problem~\eqref{inlayer1} and $\gamma\in (0,1/2]$. There exists a constant $c$ independent of $\varepsilon$ such that
\begin{equation}
\label{estiminlayer1}
\Vert \theta^{\eps}(\cdot,t)\Vert_{L^2(0,1)}  \leq c e^{-\frac{\eps^\gamma}{\eps}}
+ c\eps^{\frac{1}{2}}
e^{-\frac{M^2}{2\eps^\gamma}t} \quad \forall t\in [0,T].
\end{equation}
\end{lemma}

\noindent{\sc Proof.} i) We first check that the initial data $\theta^{\eps}(\cdot,0)$ is given by 
\begin{equation}
\label{thetat0}
\theta^{\eps}(x,0)=-\big(y_0(1)+y_0^{(1)}(1) z\big)e^{-Mz}, \quad \forall x\in (0,1].
\end{equation}
Indeed, from \eqref{def.P0eps}--\eqref{def.P3epstilde}, we write
\begin{align*}
\theta^{\eps}(x,0) & =  \sum_{k=0}^3 \eps^{k/2} \lim_{t\to 0} P_\eps^{k/2}(x,t)-y_0(x)\\
& = \lim_{t\to 0} \left(p_\eps^0(x,t) - C^\eps_0(z,\tau,t)e^{-Mz}\right)-y_0(x)\\
& \quad + \lim_{t\to 0} \sum_{k=0}^3 \eps^{k/2} \left(p^{k/2}(x,t)+Y_\eps^{k/2}\left(z,\tau,t\right) -C_{k/2}^\eps\left(z,\tau,t\right)\right), \quad x\in (0,1).
\end{align*}
We deduce from \eqref{def.W0eps}, \eqref{def.U0eps}, \eqref{def.p0tilde} and \eqref{Cztauteps}, that
\begin{align*}
& \lim_{t\to 0} \left(p_\eps^0(x,t) - C^\eps_0(z,\tau,t)e^{-Mz}\right)-y_0(x)\\
& = \lim_{t\to 0} \left(W^0_\eps(w,t)-y_0(0) - C^\eps_0(z,\tau,t)e^{-Mz}\right) = - \lim_{t\to 0} C^\eps_0(z,\tau,t)e^{-Mz}\\
& = - \lim_{t\to 0} \left(y^0(1,t)+W^0(M\tau,t)-y^0((Mt)^\pm,t)\right) e^{-Mz}=-y_0(1)e^{-Mz}.
\end{align*}
Using the matching conditions of $W^{k/2}$ with $y^k$, $1\leq k \leq 3$, we easily verify that
$\lim_{t\to 0} p^{k/2}(x,t)=0$, for $x >0$. Moreover, \eqref{Cztauteps} and \eqref{2.Y12eps.def} lead to 
\begin{align*}
& \lim_{t\to 0} \left(Y_\eps^{1/2}\left(z,\tau,t\right)-C_{1/2}^\eps\left(z,\tau,t\right)\right)\\ 
& = \lim_{t\to 0}\left(-W^{1/2}\left(M\tau,t\right)+M\tau (y^0)_x((Mt)^\pm,t)
-z W^0_{\eps,w}(M\tau,t)\right)e^{-Mz}=0, \quad \mbox{ for } x>0,
\end{align*}
then to $P_\eps^{1/2}(x,0)=0$ for all  $x\in (0,1]$. Similarly, according to \eqref{Cztauteps} and \eqref{2.Y1eps.def} we have
\begin{equation*} 
\lim_{t\to 0} \left(Y_\eps^{1}\left(z,\tau,t\right)-C_{1}^\eps\left(z,\tau,t\right)\right)
= \lim_{t\to 0}\left(-A+Bz-C^\eps\frac{z^2}{2}\right) e^{-Mz}, 
\end{equation*}
with $A$, $B$ defined by $(\ref{ABC.Y1.def})_1$ and $(\ref{ABC.Y1.def})_2$, respectively, and 
$C^\eps=W_{\eps,ww}^0(M\tau,t)$. It is easily seen that
\begin{align*}
& \lim_{t\to 0} A=0,\quad \lim_{t\to 0} C^\eps=\lim_{t\to 0} W_{\eps,ww}^0\left(M\tau,t\right)=0,\\
& \lim_{t\to 0} B=
\lim_{t\to 0}\left(-y_x^0(1,t)-W_w^{1/2}(M\tau,t)+y_x^0((Mt)^\pm,t)\right)=-y_0'(1),
\end{align*}
then $\lim_{t\to 0} \left(Y^{1}\left(z,\tau,t\right)-C^\eps_{1}\left(z,\tau,t\right)\right)
=-y_0^{(1)}(1) ze^{-Mz}$ and
\begin{equation*}
P_\eps^1(x,0)=-y_0^{(1)}(1) ze^{-Mz}, \mbox{ for } x\in (0,1].
\end{equation*}
Eventually, according to \eqref{Cztauteps} and \eqref{2.Y32eps.def} we have
\begin{equation*} 
\lim_{t\to 0} \left(Y_\eps^{3/2}\left(z,\tau,t\right)-C^\eps_{3/2}\left(z,\tau,t\right)\right)
= \lim_{t\to 0}\left(-\widetilde{A}+\widetilde{B}z-\widetilde{C}\frac{z^2}{2}+\widetilde{D}^\eps\frac{z^3}{6}\right) e^{-Mz}, 
\end{equation*}
with $\widetilde{A}$, $\widetilde{B}$, $\widetilde{C}$, defined by $(\ref{ABC.Y32.def})_1$, $(\ref{ABC.Y32.def})_2$, $(\ref{ABC.Y32.def})_3$, and $\tilde{D^\eps}= -W_{\eps,www}^0(M\tau,t)$. It is easily seen that
$$
\lim_{t\to 0} \widetilde{A}= \lim_{t\to 0} \widetilde{B}= \lim_{t\to 0} \widetilde{C}= \lim_{t\to 0} \widetilde{D}^\eps=0,
$$
then $P_\eps^{3/2}(x,0)=0$, for all $x\in (0,1]$. 

These computations lead to (\ref{thetat0}), showing that the initial condition of $\theta^\eps$ gets concentrated in the neighborhood of $x=1$. 

ii)
We now introduce a $C^{\infty}$ cut-off function $\mathcal{X}: \mathbb{R}\to [0,1]$ such that $\mathcal{X}(s)=0$ if $s\leq 1$ and $\mathcal{X}(s)=1$ if $s\geq 2$
and define, for $\gamma\in (0,1/2]$, the function $\mathcal{X}_{\eps}:[0,1]\to [0,1]$ by $\mathcal{X}_\varepsilon(x)=\mathcal{X}\left(\frac{1-x}{\varepsilon^\gamma}\right)$. 

We then decompose the solution $\theta^{\eps}$ of the linear system (\ref{inlayer1}) as folllows : $\theta^\eps=\theta^{\eps,1}+\theta^{\eps,2}$ with 
\begin{equation}
\label{inlayer1i}
\left\{
\begin{aligned}
& \theta^{\varepsilon,1}_t + M \theta^{\varepsilon,1}_x - \varepsilon \theta^{\varepsilon,1}_{xx} =0, & (x,t) \in Q_T,\\
& \theta^{\varepsilon,1}(0,t)=\theta^{\varepsilon,1}(1,t) =0, & t\in (0,T),\\
& \theta^{\varepsilon,1}(x,0)=\mathcal{X}_\varepsilon(x) \theta^{\varepsilon}(x,0), & x\in (0,1),
\end{aligned}
\right.
\qquad 
\left\{
\begin{aligned}
& \theta^{\varepsilon,2}_t + M \theta^{\varepsilon,2}_x - \varepsilon \theta^{\varepsilon,2}_{xx} =0, & (x,t) \in Q_T,\\
& \theta^{\varepsilon,2}(0,t)=\theta^{\varepsilon,2}(1,t) =0, & t\in (0,T),\\
& \theta^{\varepsilon,2}(x,0)=(1-\mathcal{X}_\varepsilon(x)) \theta^{\varepsilon}(x,0), & x\in (0,1).
\end{aligned}
\right.
\end{equation}
In view of the definition of $\mathcal{X}_{\eps}$, we see that $\theta^{\eps,1}(x,0)=0$ for all $x \geq 1-\varepsilon^\gamma$. Then, in view of (\ref{thetat0}), we check that there exists a constant $c_1>0$ independent of $\eps$ such that $\vert \theta^{\eps,1}(x,0) \vert \leq c_1 e^{-\frac{\eps^\gamma}{\eps}}$  for all $x\in(0,1)$.  By a maximum principle, it follows that
$$
\vert \theta^{\varepsilon,1}(x,t)\vert \leq c_1 e^{-\frac{\eps^\gamma}{\eps}}, \quad \forall (x,t) \in Q_T.
$$
Concerning $\theta^{\eps,2}(\cdot,0)$, we check that  $\theta^{\eps,2}(x,0)=0$ for all $x \leq 1-2\varepsilon^\gamma$ and that $\Vert \theta^{\eps,2}(\cdot,0) \Vert_{L^2(0,1)} \leq c_2 \eps^{\frac{1}{2}}$ for some constant $c_2>0$ independent of $\eps$.  Arguing as in the proof of \cite[Lemma 2.3]{amirat_munch}, we obtain that there is a constant $c_3$, independent of $\eps$, such that
$$
\Vert \theta^{\eps,2} (t,\cdot)\Vert_{L^2(0,1)} \leq c_3 \eps^{\frac{1}{2}}e^{-\frac{M^2}{2\eps^\gamma}t} \quad \forall t \in [0,T].
$$
From the two previous inequalities, we deduce \eqref{estiminlayer1}. Lemma \ref{prelim1} is proved.
$\hfill\Box$

\subsubsection{Gronwall estimate}

From now on, in order to shorten some equations, we shall use the following notation: 
\begin{equation}
L_{\eps}y:= y_t  - \eps y_{xx}+M y_x.  \nonumber
\end{equation}

We are now going to derive a priori estimates for the function $z^{\eps}$. Preliminary, since $z^{\eps}$ is not vanishing at $x=0$, we define 
\begin{equation}
\label{def:Zeps}
Z^{\eps}(x,t)=z^{\eps}(x,t)-f_\eps(x)z^{\eps}(0,t), \quad (x,t)\in Q_T,
\end{equation}
where $f_\eps$ is an appropriate function, that will be specified later, which belongs to $C^2([0,1],\mathbb{R}^+)$ and satisfies $f_\eps(0)=1$ and $f_\eps(1)=0$. The function $Z^{\eps}$ solves the following equation
\begin{equation}
\label{equation:Zeps}
	\left\{
		\begin{aligned}
   			& L_{\eps} Z^{\eps}= L_{\eps}\widetilde{P}^{\eps}-L_{\eps}(f_\eps(x) z^{\eps}(0,t)),				& (x,t)\in Q_T,    \\
   			& Z^{\eps}(0,t)  = Z^{\eps}(1,t)=0,		& t\in  (0,T), \\
			& Z^{\eps}(x, 0) = -z^{\eps}(0,0)f_\eps(x),						  		& x\in  (0,1),
   		\end{aligned} 
 	\right.
\end{equation}
where 
$z^\eps(0,0):=\lim_{t\to 0^+}z^{\eps}(0,t)$. 
We have the following result based on $L^1$-estimates with respect to the time variable.

\begin{lemma} 
\label{estimZ1} 
Let $Z^{\eps}$ be the function defined by \eqref{def:Zeps} and $\widetilde{P}^{\eps}$ the function defined by \eqref{def.P3epstilde}. There is a constant $c$ independent of $\eps$ such that, for each $t$ in $[0,T]$,
\begin{equation}
\label{estimateZ0}
\begin{aligned}
& \Vert Z^{\eps}(\cdot,t)\Vert_{L^2(0,1)}  + \sqrt{\eps} \Vert Z_x^{\eps} \Vert_{L^2((0,1)\times(0,t))} \\
& \leq c\biggl(\Vert L_{\eps}\widetilde{P}^{\eps}\Vert_{L^1(0,t,L^2(0,1))} +\Vert f_\eps\Vert_{L^2(0,1)} \biggl( \Vert z_{t}^{\eps}(0,\cdot)\Vert_{L^1(0,t)}+\vert z^{\eps}(0,0)\vert \biggr)\\
&  \quad +  \Vert -\eps f_\eps^{\prime\prime}+Mf_\eps^{\prime}\Vert_{L^2(0,1)}\Vert z^{\eps}(0,\cdot)\Vert_{L^1(0,t)} \biggr).
\end{aligned}
\end{equation}
\end{lemma}

This lemma is a consequence of the following version of the Gronwall lemma (see for instance \cite[Theorem 4, Chapter 1]{Dragomir}).

\begin{lemma} 
\label{vgl}
Let $a\in \mathbb{R}^+$, and $h \in L^2(0,T)$, $h \geq 0$. If $\zeta: [0,T] \to \mathbb{R}^+$ is a continuous function satisfying the inequality
$$\zeta(t) \leq a + 2 \int_0^t h(s)\sqrt{\zeta(s)}\,ds, \quad \forall t \in [0,T],$$ 
then we have
$$\zeta(t) \leq \left(\sqrt{a} + \int_0^t h(s)\,ds\right)^2, \quad \forall t \in [0,T].$$
\end{lemma}

\noindent \textsc{Proof of Lemma \ref{estimZ1}}- Multiplying equation \eqref{equation:Zeps} by $Z^{\eps}$ and integrating over $(0,1)\times (0,t)$ gives
\begin{equation*}
\begin{aligned}
& \frac{1}{2}\Vert Z^{\eps}(\cdot,t)\Vert^2_{L^2(0,1)} + \eps \int_0^t \Vert Z_x^{\eps}(\cdot,s)\Vert^2_{L^2(0,1)}\,ds \\
& = \frac{1}{2} \vert z^{\eps}(0,0) \vert^2 \Vert f_\eps \Vert^2_{L^2(0,1)} 
+ \int_0^t\!\!\!\int_0^1 \left(L_{\eps}\widetilde{P}^{\eps}(x,s) - L_{\eps}(f_\eps(x) z^{\eps}(0,s)) \right) Z^{\eps}(x,s)\, dxds.
\end{aligned}
\end{equation*}
Applying the Cauchy-Schwarz inequality we obtain
\begin{equation*}
\begin{aligned}
& \Vert Z^{\eps}(\cdot,t)\Vert^2_{L^2(0,1)} + 2 \eps \int_0^t \Vert Z_x^{\eps}(\cdot,s)\Vert^2_{L^2(0,1)}\,ds \leq \vert z^{\eps}(0,0) \vert^2 \Vert f_\eps \Vert^2_{L^2(0,1)}  \\
& 
+ 2\int_0^t \left(\int_0^1 \left(L_{\eps}\widetilde{P}^{\eps}(x,s) - L_{\eps}(f_\eps(x) z^{\eps}(0,s))\right)^2\, dx\right)^{1/2}
\left(\int_0^1 Z^\eps (x,s)^2\, dx\right)^{1/2}ds.
\end{aligned}
\end{equation*}
Applying Lemma \ref{vgl} with
\begin{align*}
& \zeta(t)=\Vert Z^{\eps}(\cdot,t)\Vert^2_{L^2(0,1)} + 2 \eps \int_0^t \Vert Z_x^{\eps}(\cdot,s)\Vert^2_{L^2(0,1)}\,ds, \\
& a=\vert z^{\eps}(0,0) \vert^2 \Vert f_\eps \Vert^2_{L^2(0,1)},\\
& h(s)=\left(\int_0^1 \left(L_{\eps}\widetilde{P}^\eps (x,s) - L_{\eps}(f_\eps(x) z^{\eps}(0,s))\right)^2\, dx\right)^{1/2},
\end{align*}
yields
\begin{equation*}
\begin{aligned}
& \Vert Z^{\eps}(\cdot,t)\Vert^2_{L^2(0,1)} + 2 \eps \int_0^t \Vert Z_x^{\eps}(\cdot,s)\Vert^2_{L^2(0,1)}\,ds \\
& \leq \left(\vert z^{\eps}(0,0) \vert \Vert f_\eps \Vert_{L^2(0,1)} 
+ \int_0^t \left(\int_0^1 \left(L_{\eps}\widetilde{P}^{\eps} (x,s) - L_{\eps}(f_\eps(x) z^{\eps}(0,s))\right)^2\, dx\right)^{1/2}ds\right)^2.
\end{aligned}
\end{equation*}
The previous inequality together with the equality
$$
L_{\eps}(f_\eps(x) z^{\eps}(0,t))
= f_\eps(x) z_t^{\eps}(0,t) + \left(M f^\prime_\eps(x)-\eps f^{\prime\prime}_\eps(x)\right) z^{\eps}(0,t),
$$
allow to deduce readily \eqref{estimateZ0}.
\hfill$\Box$

\subsubsection{Estimate of $\Vert z^\eps(0,\cdot)\Vert_{L^1(0,t)}$, $\Vert z^\eps_t(0,\cdot)\Vert_{L^1(0,t)}$ and $\vert z^{\eps}(0,0)\vert$}
We will now estimate each term of the right-hand side of (\ref{estimateZ0}). In the sequel, $c, c_1, c_2, \cdots$, will stand for generic constants that do not depend on $\eps$. For convenience, we recall here some notations that will be used in the sequel:
\begin{align*}
& c^+=  y_0(0), \quad c^-=v(0),\\
& d^+=  y_0^{(1)}(0), \quad d^-=-\frac{1}{M}v^{(1)}(0),\\
& e^+=y_0^{(2)}(0), \quad e^-=\frac{v^{(2)}(0)}{M^2},\\
& f^-=\frac{v^{(2)}(0)}{M^3}, \quad h^{+}=\frac{y_0^{(3)}(0)}{6}, \quad h^{-}=-\frac{v^{(3)}(0)}{6M^3}.
\end{align*} 

\begin{lemma}
\label{z(0,t)}
There is a constant $c$ independent of $\eps$ such that
\begin{equation}
\label{z(0)}
\Vert z^{\eps}(0,\cdot)\Vert_{L^1(0,t)}\leq c \eps^2, \quad \forall t>0.
\end{equation}
\end{lemma}
\textsc{Proof}- We estimate the $L^1$-norm of each term in the expansion of the function $z^{\eps,0}(0,\cdot)$. We use several time the fact that 
\begin{equation}
\Vert t^{n/2} e^{-\frac{M^2t}{4\eps}}\Vert_{L^1(0,s)}= \mathcal{O}(\eps^{1+n/2}), \quad \forall n\in \mathbb{Z}, \quad \forall s\geq 0. \label{estimatetn}
\end{equation}
Let $w_0(t)=-\frac{Mt}{\sqrt{\eps}}$. 
%
%
In view of \eqref{def.P3epstilde} and (\ref{defzeps}), $z^{\eps}(0,t)$ may be written in the form:
\begin{equation}
\label{zeps_x0}
\begin{aligned}
z^{\eps}(0,t)=&\biggl(W_\eps^0(w_0(t),t)-y^0((Mt)^-,t) \biggr)+\sqrt{\eps}\biggl(W^{1/2}(w_0(t),t)-y_x^0((Mt)^-,t)w_0(t)\biggr)\\
&+\eps\biggl(W^1(w_0(t),t)-\biggl(\frac{w_0^2}{2}y_{xx}^0((Mt)^-,t)+y^1((Mt)^-,t)\biggr)\biggr)\\
& + \eps^{3/2}\biggl(W^{3/2}(w_0(t),t)-\biggl(\frac{w^3}{3 !} y^0_{xxx}((Mt)^-,t) +  w_0(t) y_x^1( (Mt)^-,t)\biggr)\biggr)+ \mathcal{O}(e^{-\frac{M}{\eps}})\\
=& \left(p_\eps^{0}(0,t)-v(t)\right) + \sum_{k=1}^3 \eps^{k/2} p^{k/2}(0,t) + \mathcal{O}(e^{-\frac{M}{\eps}}),
\end{aligned}
\end{equation}
where $p_\eps^{0}$ is defined by \eqref{def.p0tilde}, and $p^{k/2}$, $k=1,2,3$, are defined by \eqref{def.p12}, \eqref{def.p1}, \eqref{def.p32}, respectively.
The term $\mathcal{O}(e^{-\frac{M}{\eps}})$ gathers the contributions at $x=0$ of the functions $Y^k$, introduced to describe the boundary layer at $x=1$. Far from $x=1$, these contributions are negligible. 

Observe first that from (\ref{propertyW0eps}), the first term in the right-hand side of \eqref{zeps_x0} vanishes, that is 
\begin{equation}
\label{estzeps1}
p_\eps^{0}(0,t)-v(t)=0, \quad \forall t\geq 0. 
\end{equation}
As regards the second term we have 
$$
\sqrt{\eps}p^{1/2}(0,t)=\frac{d^+-d^-}{2}(-Mt)\,erfc\left(\frac{M\sqrt{t}}{2\sqrt{\eps}}\right)+(d^+-d^-)\frac{\sqrt{\eps}\sqrt{t}}{\sqrt{\pi}}e^{-\frac{M^2t}{4\eps}}.
$$
Writing that 
$erfc(\frac{M\sqrt{t}}{2\sqrt{\eps}})=1+erf(-\frac{M\sqrt{t}}{2\sqrt{\eps}}) \leq e^{-\frac{M^2t}{4\eps}}, \; \forall t\geq 0$ using (\ref{estimate_chu}), we arrive at
$$
\sqrt{\eps} \vert p^{1/2}(0,t) \vert \leq c \vert d^+-d^- \vert \left( t + \sqrt{t}\sqrt{\eps}\right)e^{-\frac{M^2t}{4\eps}}.
$$
leading, in view of (\ref{estimatetn}) with $n=1$ and $n=2$, to
\begin{equation}
\label{estzeps2}
\sqrt{\eps} \Vert p^{1/2}(0,\cdot) \Vert_{L^1(0,t)}\leq c \vert d^+-d^- \vert  \eps^2, \quad \forall t>0.
\end{equation}
For the third term, we have
$$
\begin{aligned}
p^{1}(0,t)
=\biggl(\frac{w_0(t)^2}{2}+t\biggr)\, erfc\left(\frac{M\sqrt{t}}{2\sqrt{\eps}}\right)\frac{e^+-e^-}{2} + \frac{e^+-e^-}{2} w_0(t)\sqrt{\frac{t}{\pi}}e^{-\frac{M^2t}{4\eps}},
\end{aligned}
$$
leading to the estimate 
$$
\begin{aligned}
\eps \vert p^{1}(0,t)\vert
\leq c\vert e^+-e^-\vert \left(M^2t^2+\eps t+M\sqrt{\eps}t^{3/2} \right) e^{-\frac{M^2t}{4\eps}},
\end{aligned}
$$
then to the estimate
\begin{equation}
\label{estzeps3}
\eps \Vert p^{1}(0,\cdot) \Vert_{L^1(0,t)} \leq c\vert e^+-e^-\vert \eps^3, \quad \forall t>0.
\end{equation}
Eventually, for the fourth term, the equality
$$
\begin{aligned}
p^{3/2}(0,t) = &\biggl(\frac{w_0(t)^3}{2}+3tw_0(t)\biggr)\,erfc\left(\frac{M\sqrt{t}}{2\sqrt{\eps}}\right)(h^+-h^-)      \\
& + (h^+-h^-)(4t+w_0(t)^2)\sqrt{\frac{t}{\pi}} e^{-\frac{w_0(t)^2}{4t}}
-f^-\sqrt{\frac{t}{\pi}}  e^{-\frac{w_0(t)^2}{4t}}  - \frac{f^-}{2}w_0(t)\,erfc\left(\frac{M\sqrt{t}}{2\sqrt{\eps}}\right),
 \end{aligned}
$$
leads to 
%
$$
\begin{aligned}
&\eps^{3/2} \vert p^{3/2}(0,t) \vert
\leq c\biggl[\biggl((4t\eps^{3/2}+\eps^{1/2}M^2t^2)\sqrt{t}+\biggl(M^3t^3+3\eps Mt^2\biggr)\biggr)\vert h^+-h^-\vert  + \vert f^-\vert \biggl(\eps Mt+\eps^{3/2}\sqrt{t} \biggr)\biggr]e^{-\frac{M^2t}{4\eps}}, \\
\end{aligned}
$$
then to the estimate
\begin{equation}
\label{estzeps4}
\eps^{3/2} \Vert p^{3/2}(0,\cdot) \Vert_{L^1(0,t)}\leq c\left(\vert h^+-h^-\vert\eps^4+\vert f^-\vert \eps^3\right), \quad \forall t>0.
\end{equation}
Collecting estimates \eqref{estzeps1}--\eqref{estzeps4} we deduce from \eqref{zeps_x0}, the estimate
\eqref{z(0)}.
\hfill$\Box$

\begin{lemma}
\label{z00}
For each $\eps>0$, 
\begin{equation*}
\label{z(00)}
z^{\eps}(0,0):= \lim_{t\to 0^+} z^{\eps}(0,t) =-\biggl(y_0(1)+\frac{y_0^{(1)}(1)}{\eps}\biggr)e^{-\frac{M}{\eps}}.
\end{equation*}
\end{lemma}
\textsc{Proof-} We deduce from \eqref{def.Pk2eps}, \eqref{def.P3epstilde} and \eqref{defzeps} that
\begin{equation*}
z^{\eps}(0,t)=\widetilde{P}^{\eps}(0,t)-v(t)=\sum_{k=0}^3 \eps^{k/2} P_\eps^{k/2}(0,t)-v(t), \quad t\in (0,T). 
\end{equation*}
Due to (\ref{propertyW0eps}) we have
\begin{align*}
P_\eps^{0}(0,t)-v(t)&=-C^\eps_0\left(\frac{1}{\eps},\tau,t\right)e^{-\frac{M}{\eps}} \\
&=-\left(y_0(1-Mt)+W_\eps^0(M\tau,t)-y_0(0)\right)e^{-\frac{M}{\eps}}.
\end{align*}
Since $\lim_{t\to 0} W_\eps^0(M\tau,t)=\lim_{t\to 0} W^0(M\tau,t)=y_0(0)$, we conclude that
$$
\lim_{t\to 0} P_\eps^{0}(0,t)-v(t)=-y_0(1)e^{-\frac{M}{\eps}}.
$$
We easily verify that $\lim_{t\to 0} p^{k/2}(0,t)=0$. Morever, from \eqref{Cztauteps} and \eqref{2.Y12eps.def}
\begin{align*}
& \lim_{t\to 0} \left(Y_\eps^{1/2}\left(\frac{1}{\eps},\tau,t\right)-C_{1/2}^\eps\left(\frac{1}{\eps},\tau,t\right)\right)\\ 
& = \lim_{t\to 0}\left(-W^{1/2}\left(M\tau,t\right)+M\tau (y^0)_x((Mt)^\pm,t)
-\frac{1}{\eps} W^0_{\eps,w}(M\tau,t)\right)e^{-\frac{M}{\eps}}=0,
\end{align*}
then $\lim_{t\to 0} P_\eps^{1/2}(0,t)=0$. Similarly, according to \eqref{Cztauteps} and \eqref{2.Y1eps.def} we have
\begin{equation*} 
\lim_{t\to 0} \left(Y_\eps^{1}\left(z,\tau,t\right)-C_{1}^\eps\left(z,\tau,t\right)\right)
= \lim_{t\to 0}\left(-A+Bz-C\frac{z^2}{2}\right) e^{-Mz}, 
\end{equation*}
with $A$, $B$ defined by $(\ref{ABC.Y1.def})_1$ and $(\ref{ABC.Y1.def})_2$, respectively, and 
$C^\eps=W_{\eps,ww}^0(M\tau,t)$. It is easily seen that
\begin{align*}
& \lim_{t\to 0} A=0,\quad  \lim_{t\to 0} C^\eps=\lim_{t\to 0} W_{\eps,ww}^0\left(M\tau,t\right)=0,\\
& \lim_{t\to 0} B=
\lim_{t\to 0}\left(-y_x^0(1,t)-W_w^{1/2}(M\tau,t)+y_x^0((Mt)^\pm,t)\right)=-y_0^{(1)}(1),
\end{align*}
then $\lim_{t\to 0} \left(Y^{1}\left(z,\tau,t\right)-C_{1}\left(z,\tau,t\right)\right)
=-y_0^{(1)}(1) ze^{-Mz}$ and
\begin{equation*}
\label{P1t0}
\lim_{t\to 0} P_\eps^{1}(0,t)=0=-\frac{y_0^{(1)}(1)}{\eps} e^{-\frac{M}{\eps}}.
\end{equation*}
Eventually, according to \eqref{Cztauteps} and \eqref{2.Y32eps.def} we have
\begin{equation*} 
\lim_{t\to 0} \left(Y_\eps^{3/2}\left(z,\tau,t\right)-C^\eps_{3/2}\left(z,\tau,t\right)\right)
= \lim_{t\to 0}\left(-\widetilde{A}+\widetilde{B}z-\widetilde{C}\frac{z^2}{2}+\widetilde{D}^\eps\frac{z^3}{6}\right) e^{-Mz}, 
\end{equation*}
with $\widetilde{A}$, $\widetilde{B}$, $\widetilde{C}$, defined by $(\ref{ABC.Y32.def})_1$, $(\ref{ABC.Y32.def})_2$, $(\ref{ABC.Y32.def})_3$, and $\tilde{D^\eps}= -W_{\eps,www}^0(M\tau,t)$. It is easily seen that
$$
\lim_{t\to 0} \widetilde{A}= \lim_{t\to 0} \widetilde{B}= \lim_{t\to 0} \widetilde{C}= \lim_{t\to 0} \widetilde{D}^\eps=0,
$$
then $\lim_{t\to 0} P_\eps^{3/2}(0,t)=0$. We conclude that 
$\lim_{t\to 0^+} z^{\eps}(0,t) =-y_0(1)e^{-\frac{M}{\eps}}-\frac{y_0^{(1)}(1)}{\eps}e^{-\frac{M}{\eps}}$.
\hfill$\Box$

\begin{lemma}
\label{z_t(0)}
There is a constant $c$ independent of $\eps$ such that
\begin{equation}
\label{est:z_t(0,t)}
\Vert z_t^{\eps}(0,\cdot)\Vert_{L^1(0,t)}\leq c\,\eps, \quad \forall t>0.
\end{equation}
\end{lemma}
\textsc{Proof-} We have seen that the function $z^{\eps}(0,t)$ may be written in the form (\ref{zeps_x0}). Differentiating $z^{\eps}(0,t)$ with respect to $t$ and using the explicit form of the functions $W^k$, one can see that $z_t^{\eps}(0,t)$ may be written in the form
\begin{equation}
\label{z_t(0,t)}
z_t^{\eps}(0,t)
=\partial_t\left(p_\eps^{0}(0,t)-v(t)\right) + \sum_{k=1}^3 \eps^{k/2} p_t^{k/2}(0,t) + \mathcal{O}(e^{-\frac{M}{\eps}}).
\end{equation}
Still in view of (\ref{propertyW0eps}), 
\begin{equation}
\label{z_t(0)1}
\partial_t\left(p_\eps^{0}(0,t)-v(t)\right)=0, \quad \forall t>0,
\end{equation} 
so that the first term in the right-hand side of \eqref{z_t(0,t)} vanishes. Explicitly, for all $t>0$, 
$$
\begin{aligned}
&\sqrt{\eps} p_t^{1/2}(0,t)=\frac{d^--d^+}{2}\biggl(M\, erfc\left(\frac{M\sqrt{t}}{2\sqrt{\eps}}\right) -\frac{\sqrt{\eps}}{\sqrt{\pi}\sqrt{t}}e^{-\frac{M^2t}{4\eps}}\biggr).
\end{aligned}
$$
Therefore, for all $t>0$, 
$$
\sqrt{\eps} \left\vert p_t^{1/2}(0,t) \right\vert \leq c\vert d^--d^+\vert \left(1+\frac{\sqrt{\eps}}{\sqrt{\pi}\sqrt{t}}\right)e^{-\frac{M^2t}{4\eps}},
$$
leading to the estimate
\begin{equation}
\label{z_t(0)2}
\sqrt{\eps} \left\Vert p_t^{1/2} (0,\cdot) \right\Vert_{L^1(0,t)} \leq c\vert d^--d^+\vert \eps, 
\end{equation}
using (\ref{estimatetn}) with $n=-1$. As regard the third term, using the equality
$$
\eps p_t^{1}(0,t)=\frac{e^--e^+}{2}\biggl((M^2t+\eps)\,erfc\biggl(\frac{M\sqrt{t}}{2\sqrt{\eps}}\biggr) -\frac{2M\sqrt{\eps}\sqrt{t}}{\sqrt{\pi}}e^{-\frac{M^2t}{4\eps}}\biggr),
$$
yields
$$
\eps \bigg\vert p_t^{1}(0,t) \bigg\vert 
\leq c\vert e^--e^+\vert\biggl((M^2t+\eps)+\frac{2M\sqrt{\eps}\sqrt{t}}{\sqrt{\pi}}\biggr)e^{-\frac{M^2t}{4\eps}},
$$
then 
\begin{equation}
\label{z_t(0)3}
\eps \left\Vert p_t^{1}(0,\cdot) \right\Vert_{L^1(0,t)} \leq c\vert e^--e^+\vert \eps^2.
\end{equation}
The fourth and last term is given by 
$$
\begin{aligned}
& \eps^{3/2} p_t^{3/2}(0,t) \\
&=-\frac{3(h^+-h^-)}{\eps^4 \sqrt{\pi t}}\biggl[\sqrt{\pi}(M^3t^{5/2}\eps^{5/2}+4Mt^{3/2}\eps^{7/2})\,erfc\biggl(\frac{M\sqrt{t}}{2\sqrt{\eps}}\biggl)-e^{-\frac{M^2 t}{4\eps}}(4\eps^4t+2t^2M^2\eps^3)\biggr]\\
&\hspace*{0.5cm}+\frac{f^-}{2\eps^4 \sqrt{\pi t}}\biggl[M\eps^{7/2}\sqrt{\pi t}\,erfc\biggl(\frac{M\sqrt{t}}{2\sqrt{\eps}}\biggl)- \eps^4 e^{-\frac{M^2 t}{4\eps}}\biggr],
\end{aligned}
$$
leading to the estimate 
\begin{equation}
\label{z_t(0)4}
\eps^{3/2}\left\Vert p_t^{3/2}(0,\cdot) \right\Vert_{L^1(0,t)}
\leq c (\vert h^+-h^-\vert \eps^3 + \vert f^-\vert\eps^2).
\end{equation}
Collecting estimates \eqref{z_t(0)1}--\eqref{z_t(0)4}, we deduce from \eqref{z_t(0,t)} the estimate \eqref{est:z_t(0,t)}.
\hfill$\Box$

\subsubsection{Estimate of $\Vert L_{\eps}\widetilde{P}^{\eps}\Vert_{L^1(0,t),L^2(0,1))}$}

In order to estimate $\Vert L_{\eps}\widetilde{P}^{\eps}\Vert_{L^1(0,t),L^2(0,1))}$, we will use the following lemma. 

\begin{lemma}
\label{lemidentities}
The following identities holds:
\begin{equation}
\label{identities}
\begin{aligned}
& L_{\eps}(W^{k/2})=0, \quad 0\leq k\leq 3,\\
& L_{\eps}(e^{-Mz})=0, \quad L_{\eps}(z e^{-Mz})=\frac{M}{\eps}e^{-Mz}, \quad L_{\eps}(z^2 e^{-Mz})=-\frac{2}{\eps}(1-Mz)e^{-Mz}, \\
& L_{\eps}(z^3 e^{-Mz})=-\frac{3}{\eps}(2-Mz)ze^{-Mz} ,  \\
& L_{\eps}(w)=0,\quad L_{\eps}(w^2)=-2, \quad L_{\eps}(w^3)=-6w, \quad L_{\eps}(w^4)=-12 w^2 ,\\
& L_{\eps}(\tau)=-\frac{1}{\sqrt{\eps}}, \quad L_{\eps}(\tau^2)=-\frac{2\tau}{\sqrt{\eps}}, \quad L_{\eps}(\tau^3)= -\frac{3\tau^2}{\sqrt{\eps}}.
\end{aligned}
\end{equation}
\end{lemma}
This lemma is used to obtain the following intermediate result. 
\begin{lemma}
Let $\widetilde{P}^{\eps}$ be the function defined by \eqref{def.P3epstilde}. Assume \eqref{assumption}. Then, 
\begin{equation}
\label{LPeps}
\begin{aligned}
L_\eps(\widetilde{P}^{\eps})= & -\eps^2 y_{xx}^1(x,t) -\eps \biggl(y_t^1(1,t)+W_t^1(M\tau,t)-\partial_t(y^1((Mt)^\pm,t))\biggr)e^{-Mz} \\
& - \eps \biggl(y_{xt}^0(1,t)+W_{wt}^{1/2}(M\tau,t)\biggr)ze^{-Mz}-\eps W_{\eps,wwt}^0(M\tau,t)\frac{z^2}{2}e^{-Mz}\\
& +\eps^{3/2} \left(-\widetilde{A}_t+\widetilde{B}_t z -\widetilde{C}_t \frac{z^2}{2}+\widetilde{D}^\eps_t \frac{z^3}{6}\right)e^{-Mz},
\end{aligned}
\end{equation}
with 
\begin{equation}
\label{ABC.Y32.defb}
\left\{
\begin{aligned}
&\widetilde{A}_t=-\frac{M}{\sqrt{\eps}}W^{3/2}_{w}(M\tau,t)+W^{3/2}_{t}(M\tau,t)+ \frac{M(M\tau)^2}{2\sqrt{\eps}} y^0_{xxx}((Mt)^\pm,t) \\
&\hspace*{2cm}+ \frac{M}{\sqrt{\eps}} y_x^1( (Mt)^\pm,t)-M\tau y^0_{xxx}((Mt)^\pm,t), \\
&\widetilde{B}_t= \frac{M}{\sqrt{\eps}}W_{ww}^1(M\tau,t)- W_{wt}^1(M\tau,t) -\frac{M}{\sqrt{\eps}}y_{xx}^0((Mt)^\pm,t), \\
&\widetilde{C}_t= -\frac{M}{\sqrt{\eps}}W_{www}^{1/2}(M\tau,t)+W_{wwt}^{1/2}(M\tau,t), \\
&\widetilde{D}^\eps_t=  \frac{M}{\sqrt{\eps}}W_{\eps,wwww}^0(M\tau,t)-W_{\eps,wwwt}^0(M\tau,t).
\end{aligned}
\right.
\end{equation}
\end{lemma}
\textsc{Proof.}
We have from \eqref{def.P0eps}
\begin{equation*}
\begin{aligned}
P_\eps^0(x,t) & = p_\eps^0(x,t)+Y_\eps^0(z,\tau,t) -C^\eps_0(z,\tau,t)=p_\eps^0(x,t)-C^\eps_0(z,\tau,t)e^{-Mz}\\
& = p_\eps^0(x,t)
-\left(y^0(1,t)+W_\eps^0(M\tau,t)-y^0((Mt)^\pm,t)\right) e^{-Mz},
\end{aligned}
\end{equation*}
where $p_\eps^0(x,t)$ is defined by \eqref{def.p0tilde}. We have from \eqref{def.Pk2eps}
\begin{equation*}
\begin{aligned}
P_\eps^{1/2}(x,t)& =p^{1/2}(x,t)+Y_\eps^{1/2}(z,\tau,t) -C^\eps_{1/2}(z,\tau,t)\\
& =p^{1/2}(x,t)-\left(C^\eps_{1/2}(0,\tau,t)+W_{\eps,w}^0(M\tau,t)z\right) e^{-Mz},
\end{aligned}
\end{equation*}
where $p^{1/2}$ is defined by \eqref{def.p12}, $C^\eps_{1/2}(z,\tau,t)$ by \eqref{Cztauteps}, and $Y_\eps^{1/2}$ by \eqref{2.Y12eps.def}.
Then
\begin{equation*}
\begin{aligned}
P_\eps^0(x,t)+\sqrt{\eps} P_\eps^{1/2}(x,t)=&p_\eps^0(x,t)-\left(y^0(1,t)+W_\eps^0(M\tau,t)-y^0((Mt)^\pm,t)\right) e^{-Mz}\\
& + \sqrt{\eps} \left(p^{1/2}(x,t)-\left(C^\eps_{1/2}(0,\tau,t)+W_{\eps,w}^0(M\tau,t)z\right) e^{-Mz}\right).
\end{aligned}
\end{equation*}
Let us note here that the function $P_\eps^0+\sqrt{\eps}P_\eps^{1/2}$ belongs to $C^1([0,1]\times(0,T])$, and that $(P_\eps^0 - W_\eps^0)+\sqrt{\eps}(P_\eps^{1/2}-W^{1/2})$ belongs to $C^1(\overline{Q_T})$. Then, thanks to assumption \eqref{assumption}, $(P_\eps^0 - W_\eps^0)+\sqrt{\eps}(P_\eps^{1/2}-W^{1/2})$ belongs to 
$H^2(Q_T)$. Moreover, $L_\eps(W_\eps^0)=L^\eps(W^1)=0$. It results that, when calculating $L_\eps(P_\eps^0+\sqrt{\eps}P_\eps^{1/2})$, it suffices to perform the calculation in $\Omega^+ \cup \Omega^-$, where
$\Omega^+=\{(x,t)\in Q_T: x>Mt\}$ and $\Omega^-=\{(x,t)\in Q_T: x<Mt\}$.
A staightforward calculation then gives
\begin{equation}
\label{LP0P12}
\begin{aligned}
L_{\eps}\left(P_\eps^0+\sqrt{\eps}P_\eps^{1/2}\right)=&-\eps y_{xx}^0(x,t)-\biggl(y_t^0(1,t)-\frac{M}{\sqrt{\eps}}W_{\eps,w}^0(M\tau,t)+W_{\eps,t}^0(M\tau,t) \biggr)e^{-Mz}\\
& -\biggl(-M W_{w}^{1/2}(M\tau,t)+\sqrt{\eps}W_t^{1/2}(M\tau,t)+M y_x^0((Mt)^{\pm},t) \biggr)e^{-Mz}\\
& -\frac{M}{\sqrt{\eps}} W_{\eps,w}^0(M\tau,t) e^{-Mz} + \biggl(MW_{\eps,ww}^0(M\tau,t)-\sqrt{\eps}W_{\eps,wt}^0(M\tau,t) \biggr) ze^{-Mz}\\
=&-\eps y_{xx}^0(x,t)-\biggl(y_t^0(1,t)+W_{\eps,t}^0(M\tau,t) \biggr) e^{-Mz}\\
& -\biggl(-M W_{w}^{1/2}(M\tau,t)+\sqrt{\eps}W_t^{1/2}(M\tau,t)+M y_x^0((Mt)^{\pm},t) \biggr) e^{-Mz}\\
&  + \biggl(MW_{\eps,ww}^0(M\tau,t)-\sqrt{\eps}W_{\eps,wt}^0(M\tau,t) \biggr)ze^{-Mz}.
\end{aligned}
\end{equation}

In view of \eqref{def.p1}, \eqref{Cztauteps} and \eqref{2.Y1eps.def},
$$
P_\eps^1(x,t)=p^1(x,t)+Y_\eps^1(z,\tau,t)-C^\eps_1(z,\tau,t)=p^1(x,t)+e^{-Mz}\biggl(-A+Bz-C^\eps\frac{z^2}{2}\biggr),
$$
with 
$$
p^1(x,t)=y^1(x,t)+W^1(w,t)-\biggl(\frac{w^2}{2}y_{xx}^0((Mt)^\pm,t)+y^1((Mt)^\pm,t)\biggr),
$$
and $A$, $B$ defined by $(\ref{ABC.Y1.def})_1$ and $(\ref{ABC.Y1.def})_2$, respectively, and 
$C^\eps=W_{\eps,ww}^0(M\tau,t)$.

The function $P_\eps^{3/2}$ is given by
$$
P_\eps^{3/2}(x,t)=p^{3/2}(x,t)+Y_\eps^{3/2}(z,\tau,t)-C^\eps_{3/2}(z,\tau,t)=p^{3/2}(x,t)
+\biggl(-\widetilde{A}+\widetilde{B}z-\widetilde{C}\frac{z^2}{2}+\widetilde{D}^\eps\frac{z^3}{6}\biggr)e^{-Mz},
$$
with
$$
p^{3/2}(x,t)=W^{3/2}(w,t)-\biggl(\frac{w^3}{3 !} y^0_{xxx}((Mt)^\pm,t) + w y_x^1( (Mt)^\pm,t)\biggr),
$$
and $\widetilde{A}$, $\widetilde{B}$, $\widetilde{C}$, defined by $(\ref{ABC.Y32.def})_1$, $(\ref{ABC.Y32.def})_2$, $(\ref{ABC.Y32.def})_3$, and $\tilde{D^\eps}= -W_{\eps,www}^0(M\tau,t).$

We have the identities, valid in $\Omega^+ \cup \Omega^-$,
\begin{equation*}
\begin{aligned}
& \partial_t\left(y^1_x((Mt)^{\pm},t))\right)=y^0_{xxx}((Mt)^{\pm},t), \\
& L_\eps\left(\frac{w^3}{3 !} y^0_{xxx}((Mt)^{\pm},t)\right)=-w y^0_{xxx}((Mt)^{\pm},t), \\
& L_\eps\left(w y^1_x((Mt)^{\pm},t)\right)=w y^0_{xxx}((Mt)^{\pm},t),
\end{aligned}
\end{equation*}
where we used Lemma~\ref{lemidentities} for the last two. Arguing as for $P_\eps^0+\sqrt{\eps}P_\eps^{1/2}$,
we have by a direct calculation

\begin{equation}
\label{LP1P32}
\begin{aligned}
L_{\eps}\left(\eps P_\eps^1+\eps^{3/2}P_\eps^{3/2}\right)= & \eps y_{xx}^0(x,t)-\eps^2 y_{xx}^1(x,t) \\
& + \left(BM+C^\eps(1-Mz)\right)e^{-Mz} + \eps \left(-A_t+B_t z -C^\eps_t \frac{z^2}{2}\right)e^{-Mz}\\
& + \eps^{1/2}\left(\widetilde{B}M+\widetilde{C}(1-Mz)-\widetilde{D}^\eps(1-\frac{Mz}{2})z\right) e^{-Mz}\\
& +\eps^{3/2} \left(-\widetilde{A}_t+\widetilde{B}_t z -\widetilde{C}_t \frac{z^2}{2}+\widetilde{D}^\eps_t \frac{z^3}{6}\right)e^{-Mz}.
\end{aligned}
\end{equation}
Adding~\eqref{LP0P12} and \eqref{LP1P32} we obtain

\begin{equation*}
\label{LPeps0}
\begin{aligned}
L_\eps(\widetilde{P}^{\eps})= & -\eps^2 y_{xx}^1(x,t) -\eps^{1/2}e^{-Mz}W_t^{1/2}(M\tau,t)-\sqrt{\eps}W_{\eps,wt}^0(M\tau,t)ze^{-Mz}\\
&  +  \eps \left(-A_t+B_t z -C^\eps_t \frac{z^2}{2}\right)e^{-Mz}\\
& + \eps^{1/2}\biggl(\widetilde{B}M+\widetilde{C}(1-Mz)-\widetilde{D}^\eps(1-\frac{Mz}{2})z\biggr) e^{-Mz}\\
& +\eps^{3/2} \left(-\widetilde{A}_t+\widetilde{B}_t z -\widetilde{C}_t \frac{z^2}{2}+\widetilde{D}^\eps_t \frac{z^3}{6}\right)e^{-Mz}.
\end{aligned}
\end{equation*}
Rearranging the terms we arrive at (\ref{LPeps}).$\hfill\Box$

\begin{lemma}
\label{lemmeLPeps}
Let $\widetilde{P}^{\eps}$ be the function defined by \eqref{def.P3epstilde}. Assume \eqref{assumption}. Then there is a constant $c$ independent of $\eps$ such that
\begin{equation}
\label{LP3epsf}
\Vert L_\eps\widetilde{P}^{\eps} \Vert_{L^1(0,T, L^2(0,1)} \leq c \,\eps^{3/2}.
\end{equation}
\end{lemma}

\noindent
\textsc{Proof.} We estimate the $L^1(0,t; L^2(0,1)$-norm of each term of the right-hand side of \eqref{LPeps}. We use notably several times that 
\begin{itemize}
\item $\Vert z(x)^n e^{-Mz(x)}\Vert_{L^2(0,1)}=\mathcal{O}(\sqrt{\eps})$ for all $n\in \mathbb{N}$, with $z=(1-x)/\eps$. 
\item $\biggl\Vert (1-Mt)^n e^{-\frac{(1-Mt)^2}{4\eps t}}\biggr\Vert_{L^1(0,T)}=\mathcal{O}(\eps^{(n+1)/2})$ for all $n\in \mathbb{N}$.
\end{itemize}
\par\noindent
{\textbf{a)}We have $\Vert -\eps^2 y_{xx}^1 \Vert_{L^2(0,t; L^2(0,1))} \leq c \eps^2$, then, using the Cauchy-Schwarz inequality, 
\begin{equation}
\label{estim001}
\Vert -\eps^2 y_{xx}^1 \Vert_{L^1(0,t; L^2(0,1))} \leq c \eps^2.
\end{equation}
\par\noindent 
{\textbf{b)} We have by a direct calculation
\begin{equation}
\label{estim002}
\Vert -\eps y^1(1,t)e^{-Mz}\Vert_{L^1(0,t,L^2(0,1))}=\eps \Vert y^1(1,t)\Vert_{L^1(0,t)}\Vert e^{-Mz}\Vert_{L^2(0,1)}\leq c \eps^{3/2}.
\end{equation}
\par\noindent 
{\textbf{c)} The estimation of $\Vert W_t^1(M\tau,t)-\partial_t(y^1((Mt)^\pm,t)\Vert_{L^1(0,t)}$ can be achieved as follows. Explicitly
$$
W_t^1(w,t)=\frac{e^+-e^-}{2}erf\biggl(\frac{w}{2\sqrt{t}}\biggr)+\frac{e^++e^-}{2},\quad \partial_t(y^1((Mt)^\pm,t)=e^{\pm},
$$ 
then
$$
W_t^1(M\tau,t)-\partial_t(y^1((Mt)^\pm,t)=
\left\{
\begin{aligned}
&\frac{1}{2}(e^--e^+)\,erfc\biggl(\frac{1-Mt}{2\sqrt{t}\sqrt{\eps}} \biggr), & Mt<1, \\
&\frac{1}{2}(e^+-e^-)\,erfc\biggl(\frac{Mt-1}{2\sqrt{t}\sqrt{\eps}} \biggr), & Mt>1,
\end{aligned}
\right.
$$
then
\begin{equation}
\label{estim003}
\Vert W_t^1(M\tau,t)-\partial_t(y^1((Mt)^\pm,t)\Vert_{L^1(0,t)}\leq \vert e^+-e^-\vert  \Vert e^{-\frac{(1-Mt)^2}{4\eps t}}\Vert_{L^1(0,t)}\leq c \vert e^+-e^-\vert \eps^{1/2}.
\end{equation}
We conclude that, 
\begin{equation}
\label{estim004}
\biggl\Vert \eps \biggl(W_t^1(M\tau,t)-\partial_t(y^1((Mt)^\pm,t))\biggr)e^{-Mz}\biggr\Vert_{L^1(0,T,L^2(0,1))}\leq c \vert e^+-e^-\vert  \eps^2.
\end{equation}
\par\noindent 
{\textbf{d)}
A direct calculation gives
\begin{equation}
\label{estim005}
\Vert -\eps y_{xt}^0(1,t) ze^{-Mz}\Vert_{L^1(0,T,L^2(0,1))} \leq c \eps^{3/2}.
\end{equation}

Let us then estimate $\Vert W_{wt}^{1/2}(M\tau,t)\Vert_{L^1(0,T)}$. Explicitly, $W_{wt}^{1/2}(w,t)=\frac{1}{4\sqrt{\pi}t^{3/2}}(d^--d^+)we^{-\frac{w^2}{t}}$, then 
$$
W_{wt}^{1/2}(M\tau,t)=\frac{d^--d^+}{4\sqrt{\eps}\sqrt{\pi}t^{3/2}}(1-Mt)e^{-\frac{(1-Mt)^2}{4\eps t}},
$$
hence $\Vert W_{wt}^{1/2}(M\tau,t)\Vert_{L^1(0,t)}\leq c \vert d^+-d^-\vert \sqrt{\eps}$. It then follows
\begin{equation}
\label{estim006}
\eps \Vert W_{wt}^{1/2}(M\tau,t) z e^{-Mz}\Vert_{L^1(0,T,L^2(0,1))} \leq c \vert d^+-d^-\vert \eps^2.
\end{equation}
\par\noindent 
{\textbf{e)}
For the estimation of $\Vert W_{\eps,wwt}^0(M\tau,t)\Vert_{L^1(0,t)}$ we write $W_{\eps,wwt}^0(M\tau,t)=W_{wwt}^0(M\tau,t)+{U}_{\eps,wwt}^0(M\tau,t)$.
Estimation of the first term is direct. We have
$$
W_{wwt}^0(w,t)=\frac{c^+-c^-}{16t^{7/2}\sqrt{\pi}}w(6t-w^2)e^{-\frac{w^2}{4t}},
$$
then
$$
W_{wwt}^0(M\tau,t)=\frac{c^+-c^-}{16t^{7/2}\sqrt{\pi}\eps^{3/2}}(1-Mt)(-6\eps t-(1-Mt)^2)e^{-\frac{(1-Mt)^2}{4\eps t}},
$$
from which we deduce that $\Vert W_{wwt}^0(M\tau,t)\Vert_{L^1(0,T)}\leq c \vert c^+-c^-\vert$, hence
\begin{equation}
\label{estim007}
\left\Vert \eps W_{wwt}^0(M\tau,t)\frac{z^2}{2}e^{-Mz}\right\Vert_{L^1(0,t,L^2(0,1))} \leq c \vert c^+-c^-\vert \eps^{3/2}.
\end{equation}
Estimate of the second requires more care. A straightforward calculation gives
\begin{equation*}
\begin{aligned}
& U_{\eps,wwt}^0(M\tau,t)=-\frac{c^+-c^-}{8} e^{\frac{M}{\eps}}\frac{1}{\eps^{2}\sqrt{\pi}t^{7/2}}\\
& \times \left(-8M^4 \sqrt{\pi}t^{7/2}\,erfc\biggl(\frac{1+Mt}{2\sqrt{t}\sqrt{\eps}}\biggr)+\eps^{1/2} e^{-\frac{(1+Mt)^2}{4t\eps}}(-11t^2M^2 +15t^3 M^3 +5tM -10\eps Mt^2 +6\eps t-1)\right).
\end{aligned}
\end{equation*}
We now use the asymptotic behavior (\ref{behavior_erf}) of the erf function to write that 
$$
-\sqrt{\pi}\,erfc\left(\frac{1+Mt}{2\sqrt{t}\sqrt{\eps}}\right)=e^{-\frac{(1+Mt)^2}{4t\eps}}\left(-\frac{2\sqrt{\eps}\sqrt{t}}{1+Mt}+\mathcal{O}(\eps^{3/2})\right),
$$
then 
\begin{equation*}
\begin{aligned}
&\left(-8M^4 \sqrt{\pi}t^{7/2}\, erfc\biggl(\frac{1+Mt}{2\sqrt{t}\sqrt{\eps}}\biggr)+\sqrt{\eps}e^{-\frac{(1+Mt)^2}{4t\eps}}(-11t^2M^2 +15t^3 M^3 +5tM -10\eps Mt^2 +6\eps t-1) \right)\\
&=e^{-\frac{(1+Mt)^2}{4t\eps}}\left(\sqrt{\eps}\left(-\frac{16M^4t^4}{1+Mt}-11t^2M^2 +15t^3 M^3 +5tM -1\right)+ \mathcal{O}(\eps^{3/2})\right)\\
&=e^{-\frac{(1+Mt)^2}{4t\eps}}\left(-\sqrt{\eps}\frac{(1-Mt)^4}{(1+Mt)}+ \mathcal{O}(\eps^{3/2})\right).
\end{aligned}
\end{equation*}
Since $e^{\frac{M}{\eps}}e^{-\frac{(1+Mt)^2}{4t\eps}}=e^{-\frac{(1-Mt)^2}{4t\eps}}$, we may write 
\begin{equation*}
U_{\eps,wwt}^0(M\tau,t)=\frac{c^+-c^-}{8}\frac{1}{\eps^{2}\sqrt{\pi}t^{7/2}}e^{-\frac{(1-Mt)^2}{4t\eps}}\left(\sqrt{\eps}\frac{(1-Mt)^4}{(1+Mt)}+ \mathcal{O}(\eps^{3/2})\right),
\end{equation*}
leading to $\Vert U_{\eps,wwt}^0(M\tau,t)\Vert_{L^1(0,T)}\leq c \vert c^+-c^-\vert$, then to
\begin{equation}
\label{estim008}
\left\Vert \eps U_{\eps,www}^0(M\tau,t)\frac{z^2}{2}e^{-Mz}\right\Vert_{L^1(0,T,L^2(0,1))} \leq c \vert c^+-c^-\vert \eps^{3/2}.
\end{equation}
Adding \eqref{estim007} and \eqref{estim008} we obtain
\begin{equation}
\label{estim009}
\left\Vert \eps W_{\eps,www}^0(M\tau,t)\frac{z^2}{2}e^{-Mz}\right\Vert_{L^1(0,T,L^2(0,1))} \leq c \vert c^+-c^-\vert \eps^{3/2}.
\end{equation}
\par\noindent 
{\textbf{f)} Let us now estimate $\Vert \widetilde{A}_t\Vert_{L^1(0,t)}$. Straightforward calculations give
$$
\widetilde{A}_t(M\tau,t)=
\left\{
\begin{aligned}
& (h^+-h^-)\biggl(\frac{3\sqrt{t}}{\eps
\sqrt{\pi}}e^{-\frac{(1-Mt)^2}{4t\eps}}(2\eps-M(1-Mt)) \biggr) \\
& + \frac{3}{2\eps^{3/2}}(h^+-h^-)\biggl(2\eps
(2Mt-1)+M(1-Mt)^2\biggr)\, erfc\biggl(\frac{1-Mt}{2\sqrt{\eps}\sqrt{t}}\biggr)\\
&
-\frac{f^-}{2}\biggl(\frac{1}{\sqrt{t}\sqrt{\pi}}e^{-\frac{(1-Mt)^2}{4t\eps}}+\frac{M}{\sqrt{\eps}}\, erfc\biggl(\frac{1-Mt}{2\sqrt{\eps}\sqrt{t}}\biggr)\biggr),
\qquad \qquad 1-Mt>0,\\\\
& (h^+-h^-)\biggl(\frac{3\sqrt{t}}{\eps
\sqrt{\pi}}e^{-\frac{(1-Mt)^2}{4t\eps}}(2\eps-M(1-Mt)) \biggr) \\
& + \frac{3}{2\eps^{3/2}}(h^+-h^-)\biggl(2\eps
(1-2Mt)+M(1-Mt)^2\biggr)\, erfc\biggl(\frac{Mt-1}{2\sqrt{\eps}\sqrt{t}}\biggr)\\
&
-\frac{f^-}{2}\biggl(\frac{1}{\sqrt{t}\sqrt{\pi}}e^{-\frac{(1-Mt)^2}{4t\eps}}-\frac{M}{\sqrt{\eps}}\,erfc\biggl(\frac{Mt-1}{2\sqrt{\eps}\sqrt{t}}\biggr)\biggr),
\qquad \qquad 1-Mt<0,\\
\end{aligned}
\right.
$$
leading to $\Vert \widetilde{A}_t(M\tau,t)\Vert_{L^1(0,T)}\leq c (\vert
h^+-h^-\vert+\vert f^-\vert)$, then to
\begin{equation}
\label{estim010}
\Vert \eps^{3/2}\widetilde{A}_t(M\tau,t)e^{-Mz}\Vert_{L^1(0,t,L^2(0,1))}\leq c
(\vert h^+-h^-\vert+\vert f^-\vert) \eps^2.
\end{equation}
\par\noindent 
{\textbf{g)} To estimate $\Vert \widetilde{B}_t\Vert_{L^1(0,t)}$ we use the equality
$$
\widetilde{B}_t(w,t)=
\left\{
\begin{aligned}
&\frac{M}{2\sqrt{\eps}} \biggl(erf\biggl(\frac{w}{2\sqrt{t}}\biggr)-1\biggr)(e^+-e^-)-\frac{1}{2}\frac{e^{-\frac{w^2}{4t}}}{\sqrt{\pi}\sqrt{t}}(e^+-e^-),  \quad w>0,\\
&\frac{M}{2\sqrt{\eps}} \biggl(erf\biggl(\frac{w}{2\sqrt{t}}\biggr)+1\biggr)(e^+-e^-)-\frac{1}{2}\frac{e^{-\frac{w^2}{4t}}}{\sqrt{\pi}\sqrt{t}}(e^+-e^-),  \quad w<0,\\
\end{aligned}
\right.
$$
leading to 
$$
\vert \widetilde{B}_t(M\tau,t)\vert \leq 
\biggl(\frac{M}{2\sqrt{\eps}} +\frac{1}{2\sqrt{\pi}\sqrt{t}}\biggr) e^{-\frac{(1-Mt)^2}{4\eps t}}\vert e^+-e^-\vert,
$$
then to $\Vert \widetilde{B}_t(M\tau,t)\Vert_{L^1(0,T)}\leq c$, and then to 
\begin{equation}
\label{estim011}
\Vert \eps^{3/2}\widetilde{B}_t(M\tau,t)ze^{-Mz}\Vert_{L^1(0,t,L^2(0,1))}\leq c \vert e^+-e^-\vert \eps^2.
\end{equation}
\par\noindent 
{\textbf{h)} To estimate $\Vert \widetilde{C}_t\Vert_{L^1(0,T)}$ we write
$$
\begin{aligned}
\widetilde{C}_t(M\tau,t)&=-\frac{M}{\sqrt{\eps}}W_{www}^{1/2}(M\tau,t)+W_{wwt}^{1/2}(M\tau,t)\\
&=-\frac{d^+-d^-}{\sqrt{\pi}t^{5/2}\eps}\biggl( (1-Mt)(1+Mt) + 2\eps t\biggr)e^{-\frac{(1-Mt)^2}{4t\eps}},
\end{aligned}
$$
then 
$\Vert \widetilde{C}_t(M\tau,t)\Vert_{L^1(0,T)}\leq c \vert d^+-d^-\vert$, hence
\begin{equation}
\label{estim012}
\Vert \eps^{3/2}\widetilde{C}_t(M\tau,t)z^2e^{-Mz}\Vert_{L^1(0,T,L^2(0,1))}\leq c \vert d^+-d^-\vert \eps^2.
\end{equation}
\par\noindent 
{\textbf{i)} It remains to estimate $\Vert \widetilde{D}_t\Vert_{L^1(0,t)}$. We have
$$
\begin{aligned}
\widetilde{D}_t(M\tau,t)&=\frac{M}{\sqrt{\eps}}W_{\eps,wwww}^0(M\tau,t)-W_{\eps,wwwt}^0(M\tau,t)\\
&=\frac{c^+-c^-}{16\sqrt{\pi}t^{9/2}\eps^2} \biggl(6\eps t\biggl[(1-Mt)^2+(1-Mt)-2t\eps\biggr]-(1-Mt)^3 \biggr)e^{-\frac{(1-Mt)^2}{4t\eps}},
\end{aligned}
$$
then $\Vert \widetilde{D}_t(M\tau,t)\Vert_{L^1(0,T)}\leq c \vert c^+-c^-\vert$, hence
\begin{equation}
\label{estim013}
\Vert \eps^{3/2}\widetilde{D}_t(M\tau,t)z^3e^{-Mz}\Vert_{L^1(0,T,L^2(0,1))}\leq c \vert c^+-c^-\vert \eps^2.
\end{equation}
Collecting estimates \eqref{estim001}--\eqref{estim013} we deduce from \eqref{LPeps} the estimate \eqref{LP3epsf}. The proof of Lemma \ref{lemmeLPeps} is complete. 
\hfill$\Box$

\subsection{End of the proof of Theorem \ref{convergence}}

    We are now in position to finish the proof of Theorem \ref{convergence}. It remains to choose the function $f_\eps\in C^2([0,1])$, satisfying $f_\eps(0)=1$ and $f_\eps(1)=0$ so as to minimize the terms in the right side of (\ref{estimateZ0}), asymptotically with respect to $\eps$. Since the terms $\Vert z_{t}^{\eps}(0,\cdot)\Vert_{L^1(0,t)}+\vert z^{\eps}(0,0)\vert$  and $\Vert z^{\eps}(0,\cdot)\Vert_{L^1(0,t)}$ are respectively of order $\eps$ and $\eps^2$, according to Lemmas \ref{z(0,t)}, \ref{z00}, and \ref{z_t(0)}, we consider the  function
$$
f_\eps(x)=(1-x) e^{-\frac{Mx}{\eps}}, \quad x\in [0,1],
$$
so that $\Vert f_\eps \Vert_{L^2(0,1)} \leq c \,\eps^{1/2}$ and  $\Vert -\eps f_\eps^{\prime\prime}+Mf_\eps^{\prime} \Vert_{L^2(0,1)} \leq c \,\eps^{-1/2}$. It follows that 
$$
\Vert f_\eps\Vert_{L^2(0,1)} \left( \Vert z_{t}^{\eps}(0,\cdot)\Vert_{L^1(0,t)}+\vert z^{\eps}(0,0)\vert \right)
\leq c \eps^{3/2},
$$
$$
\Vert -\eps f_\eps^{\prime\prime}+Mf_\eps^{\prime}\Vert_{L^2(0,1)}\Vert z^{\eps}(0,\cdot)\Vert_{L^1(0,t)}
\leq c \eps^{3/2}.
$$
Coming back to Lemma \ref{estimZ1}, using the two previous estimates and Lemma \ref{lemmeLPeps}, it follows that
$$
\Vert Z^{\eps}(\cdot,t)\Vert_{L^2(0,1)}  + \sqrt{\eps} \Vert Z_x^{\eps} \Vert_{L^2((0,1)\times(0,t))}
\leq c \eps^{3/2}, \quad \forall t\in [0,T].
$$
Now, since
$$z^{\eps}(x,t)=Z^{\eps}(x,t)+f_\eps(x)z^{\eps}(0,t),$$
we deduce that
\begin{equation}
\label{zinfty}
\Vert z^{\eps}(\cdot,t)\Vert_{L^2(0,1)}  
\leq \Vert Z^{\eps}(\cdot,t)\Vert_{L^2(0,1)} + \vert z^{\eps}(0,t) \vert \Vert f_\eps \Vert_{L^2(0,1)}
\leq c \eps^{\frac{3}{2}} + c \vert z^{\eps}(0,t)  \vert \eps^{\frac{1}{2}}, \quad \forall t\in [0,T].
\end{equation}
Then, writing
$$z^{\eps}(0,t)=z^{\eps}(0,0) + \int_0^t z_t^{\eps}(0,s)\, ds,$$
and using Lemmas \ref{z00} and \ref{z_t(0)}, it holds that
\begin{equation*}
\vert z^{\eps}(0,t) \vert \leq  \vert z^{\eps}(0,0) \vert  + \Vert z_t^{\eps}(0,\cdot) \Vert_{L^1(0,t)} \\
\leq c\left(1+\frac{1}{\eps}\right) e^{-\frac{M}{\eps}} + c\eps, \quad \forall t>0.
\end{equation*}
This estimate allows deduce from \eqref{zinfty} that
\begin{equation*}
\Vert z^{\eps}(\cdot,t)\Vert_{L^2(0,1)}  
\leq c \eps^{\frac{3}{2}} + c\eps^{\frac{1}{2}}
\left(\left(1+\frac{1}{\eps}\right) e^{-\frac{M}{\eps}} + \eps\right) \leq c \eps^{\frac{3}{2}}, \quad \forall t\in [0,T].
\end{equation*}
Then, since $\widetilde{P}^{\eps}(x,t)-y^{\eps}(x,t)=z^{\eps}(x,t)-\theta^{\eps}(x,t)$, using Lemma \ref{prelim1}, we derive the estimate
\begin{align*}
\Vert \widetilde{P}^{\eps}(\cdot,t)-y^{\eps}(\cdot,t) \Vert_{L^2(0,1)} & \leq c\, \eps^{3/2} + c e^{-\frac{\eps^\gamma}{\eps}}
+ c\, \eps^{1/2}
e^{-\frac{M^2}{2\eps^\gamma}t} \\
& \leq c\,\eps^{3/2} 
+ c\,\eps^{1/2} e^{-\frac{M^2}{2\eps^\gamma}t}\quad \forall t\in [0,T].
\end{align*}
This ends the proof of  Theorem \ref{convergence}. $\hfill\Box$

\subsection{Remarks}

The following remarks are in order. 

\begin{remark}  We have constructed the approximation $\widetilde{P}^\eps$ in two steps: first, we have used the matching asymptotic method to derive a composite approximation $P^\eps$, as a non trivial linear combination of the functions $y^k$, $W^k$, and $Y^k$, $k=0,\cdots,3$. Then, in order to get a better estimate and eliminate an artificial boundary layer propagating along the characteristic, we have defined $\widetilde{P}^\eps$ by replacing the function $W^0$ by the function $W^0_\eps=W^0+U^0_\eps$. This latter depends explicitly on $\eps$ and is associated to the integral representation of the exact solution of the advection-diffusion system defined over $\mathbb{R}^+\times \mathbb{R}^+$.  We may also proceed directly with the function $W_0^\eps$ through the method of matching asymptotic expansion. The methodology is the same but since $W^0_\eps$ depends explicitly on $\eps$, we need to expand $W^0_\eps$ when we determine the matching conditions for the inner layer (see section \ref{innerexp}) and the boundary one (see section \ref{innerexpansionx1}). The introduction of $W^0_\eps$ does not modify the matching conditions (\ref{matching_condition_w}) for the  inner layer since $W^0$ and $W^0_\eps$ an their derivatives share the same asymptotic behavior as $w \to \pm \infty$; we check that  
\begin{equation}
U^0_{\eps}(w,t)=\frac{\eps^{1/2}}{M\sqrt{t \pi}} \biggl(1-\eps^{1/2}\frac{w}{2 tM} +\eps\frac{(w^2-2t)}{4t^2M^2}+\cdots\biggr)e^{-\frac{w^2}{4t}} \label{expansionU0eps}
\end{equation}
so that $U^0_{\eps}(w,t)\to 0$ as $w\to\pm \infty$, for all $t>0$. On the other hand, this modifies the matching conditions (\ref{matching_condition_z}) since $W^0$ and $W^0_\eps$ does not share necessarily the same limit as $z\to \infty$: 
$$
\begin{aligned}
U^0_\eps(w,t)  =U^0_\eps(M\tau-\eps z,t)=& U^0_\eps(M\tau,t) - \eps z U^0_{\eps,z}(M\tau,t) + \cdots\\
=& \eps^{1/2}\frac{1}{M\sqrt{t\pi}}e^{-\frac{(M\tau)^2}{4t}} +\mathcal{O}(\eps)
\end{aligned}
$$
where we have used (\ref{expansionU0eps}) with $w=M\tau$. This implies additional terms in the definition of the functions $C_k(z,\tau,t)$, defined in (\ref{Cztaut}) (used to construct $P^\eps$) and therefore some changes in the definition of the function $C^\eps_k(z,\tau,t)$ defined in (\ref{Cztauteps}) (used to construct $\widetilde{P}^\eps$). For instance, since $U^0_\eps$ has no contribution to the power $\eps^0$, $C^{\eps}_0(z,\tau,t)=y^0(1,t)+W_\eps^0(M\tau,t)-y^0((Mt)^\pm,t)$ becomes simply $C^{\eps}_0(z,\tau,t)=y^0(1,t)+W^0(M\tau,t)-y^0((Mt)^\pm,t)$, that is $C_0(z,\tau,t)$. On the other hand, $C_{1/2}^\eps$ defined in (\ref{Cztauteps}) as follows
$$
C_{1/2}^\eps(z,\tau,t)=W^{1/2}(M\tau,t)-M\tau (y^0)_x((Mt)^\pm,t)-z W_{\eps,w}^0(M\tau,t)
$$
has to be replaced by
$$
C_{1/2}^\eps(z,\tau,t)=W^{1/2}(M\tau,t)-M\tau (y^0)_x((Mt)^\pm,t)+\frac{1}{M\sqrt{t\pi}}e^{-\frac{(M\tau)^2}{4t}}-z W_w^0(M\tau,t).
$$
Both approaches lead to the same estimate and requires a similar quantity of calculus.  
\end{remark}

\begin{remark}
Concerning the error estimate obtained when we use $W^0$ leading to the approximation $P^\eps$ defined in (\ref{def.P3eps}), we compute that 
$$
\biggl\vert W^0\biggl(-\frac{Mt}{\sqrt{\eps}},t\biggr)-y^0((Mt)^-,t)\biggr\vert=\frac{\vert y_0(0)-v(0)\vert}{2}\,erfc\biggl(\frac{M\sqrt{t}}{2\sqrt{\eps}}\biggr) \leq \vert y_0(0)-v(0)\vert e^{-\frac{M^2t}{4\eps}}, \quad \forall t\geq 0
$$
leading the estimate $\Vert W^0(-\frac{Mt}{\sqrt{\eps}},t)-y^0((Mt)^-,t)\Vert_{L^1(0,T)}=\mathcal{O}(\eps)$, while $W_\eps^0(-\frac{Mt}{\sqrt{\eps}},t)-y^0((Mt)^-,t)$ vanishes for all $t\geq 0$.
Similarly, we explicitly compute that 
$$
\biggl\Vert \partial_t\biggl(W^0\biggl(-\frac{Mt}{\sqrt{\eps}},t\biggr)-y^0((Mt)^-,t)\biggr)\biggr\Vert_{L^1(0,s)}=\vert y_0(0)-v(0)\vert \,erf\biggl(\frac{M\sqrt{s}}{2\sqrt{\eps}}\biggr)\to 1\quad\textrm{as}\quad\eps\to  0,\quad\forall s>0.
$$
It results that the analogous of $z^{\eps}$ (defined in (\ref{defzeps})  with $W^0_\eps$) associated to $W^0$ satisfies $\Vert z^{\eps}(0,\cdot)\Vert_{L^1(0,t)}=\mathcal{O}(\eps^{1/2})$ and  
$\Vert z_t^{\eps}(0,\cdot)\Vert_{L^1(0,t)}=\mathcal{O}(1)$. Gronwall estimate of Lemma \ref{estimZ1} then leads to 
$$
\Vert P^\eps-y^{\eps}\Vert_{L^\infty(0,T,L^2(0,1))}=\vert y_0(0)-v(0)\vert \mathcal{O}(\eps^{1/2})
$$
 to be compared to the rate $\mathcal{O}(\eps^{3/2})$ obtained when $W^0_\eps$ is used.  Remark however that if $y_0(0)=v(0)$ (absence of shock layer), then $W^0=W^0_\eps=0$ and the estimate is of order $\mathcal{O}(\eps^{3/2})$.

Coming back to the error estimate on $\widetilde{P}^{\eps}$ based on $W^0_\eps$, we see that the rate is driven by the $L^1$-norm of the time derivative at $x=0$, i.e. $\Vert z_t^\eps(0,\cdot)\Vert_{L^1(0,t)}$ of the order $\eps$, see Lemma \ref{z_t(0)}. Precisely, in the expansion of $z_t^\eps(0,\cdot)$, the main term comes from the function 
$$
p^{1/2}(0,t)= W^{1/2}\biggl(-\frac{Mt}{\sqrt{\eps}},t\biggr)-d^{-}\biggl(-\frac{Mt}{\sqrt{\eps}}\biggr).
$$
see (\ref{z_t(0)2}). In order to improve the estimate, we may modify the function $W^{1/2}$ and replace it by the function 
$W^{1/2}_\eps=W^{1/2} + U^{1/2}_\eps$ with 
\begin{equation}
\label{def.U12eps}
U^{1/2}_\eps(w,t)=(d^+-d^-) \sqrt{t}e^{\frac{Mw}{\sqrt{\eps}}+\frac{M^2t}{\eps}}ierfc\left(\frac{w}{2\sqrt{t}}+\frac{M\sqrt{t}}{\sqrt{\eps}}\right),
\end{equation}
and $ierfc(w):=\int_w^\infty erfc(s)ds=\frac{1}{\sqrt{\pi}}e^{-w^2}-w\, erfc(w)$, for all $w \in \mathbb{R}$. In particular, we check that $W_{\eps}^{1/2}$ satisfy 
\begin{equation}
W_\eps^{1/2}\biggl(-\frac{Mt}{\sqrt{\eps}},t\biggr)-d^- \biggl(-\frac{Mt}{\sqrt{\eps}}\biggr)=0, \quad \forall t\in (0,T]. \label{propertyW12eps}
\end{equation}
to be compared with (\ref{limit_t0W12}). Actually, the function $\widetilde{W}^{1/2}_{\eps}(x,t)=W^{1/2}_{\eps}(\frac{x-Mt}{\sqrt{\eps}},t)$ solves the equation (we refer notably to \cite{Shih07})
\begin{equation*}
	\left\{
		\begin{aligned}
   			& \widetilde{W}^{1/2}_{\eps,t} + M \widetilde{W}^{1/2}_{\eps,x} - \eps \widetilde{W}^{1/2}_{\eps,xx}=0, 				& (x,t)\in \mathbb{R}^+\times \mathbb{R}^+,    \\
   			& \widetilde{W}^{1/2}_{\eps}(0,t)  = d^- \biggl(-\frac{Mt}{\sqrt{\eps}}\biggr), 	  		& t\in  \mathbb{R}^+, \\
   			& \widetilde{W}^{1/2}_{\eps}(x, 0) = d^+ \biggl(\frac{x}{\sqrt{\eps}}\biggr), 						  		& x\in  \mathbb{R}^+.
   		\end{aligned} 
 	\right.
\end{equation*}
As a consequence, the use of the function $p_\eps^{1/2}(x,t):=W^{1/2}_\eps(w,t)-y^0( (Mt)^{\pm},t)w$ ($w=(x-Mt)/\sqrt{\eps}$) instead of $p^{1/2}$ defined in (\ref{def.p12}) allows to improve the approximation $\widetilde{P}^\eps$. 
\end{remark}

\begin{remark}
Gronwall estimate (\ref{estimateZ0}) implies that $\sqrt{\eps}\Vert Z^\eps_x\Vert_{L^2((0,1)\times (0,t)}=\mathcal{O}(\eps^{3/2})$. Writing that  
\begin{equation*}
\sqrt{\eps}\Vert (\widetilde{P}^\eps-y^{\eps})_x\Vert_{L^2(Q_t)}\leq \sqrt{\eps}\Vert Z^\eps_x\Vert_{L^2((0,1)\times (0,t)}  +\sqrt{\eps}\Vert \theta_x^\eps\Vert_{L^2(Q_t)} +\sqrt{\eps}\Vert f^\eps_x z^{\eps}(0,\cdot)\Vert_{L^2(Q_t)},
\end{equation*}
that $\Vert f^\eps_x \Vert _{L^2(0,1)}=\mathcal{O}(\eps^{-1/2})$ and  $\Vert z^\eps(0,\cdot) \Vert_{L^2(0,t)} =\mathcal{O}(\eps^{3/2})$ leading to $\Vert f^\eps_x z^{\eps}(0,\cdot)\Vert_{L^2(Q_t)} = \mathcal{O(\eps)}$, we deduce that there exists a constant $c>0$ such that 
\begin{equation*}
\Vert (\widetilde{P}^\eps-y^{\eps})_x\Vert_{L^2(Q_T)}\leq c\,\eps + \Vert \theta_x^\eps\Vert_{L^2(Q_T)}
\end{equation*}
where $\theta$ is the initial layer corrector defined in (\ref{inlayer1}), decomposed as follows $\theta^\eps=\theta^{\eps,1}+\theta^{\eps,2}$, see (\ref{inlayer1i}).
In view of the structure of $\theta^{\eps,1}(\cdot,0)$, we easily show that $\Vert \theta^{\eps,1}_x\Vert_{L^2(Q_T)}\leq c \eps^{-1/2}e^{-\eps^{\gamma-1}}$ for all $\gamma\in (0,1/2]$.
On the other hand, energy estimate for $\theta^{\eps,2}$, leads, in view of (\ref{thetat0}), to  
$$
\begin{aligned}
\Vert \theta_x^{\eps,2}\Vert_{L^2(Q_T)}& \leq \eps^{-1/2}\Vert \theta^{\eps,2}(\cdot,0)\Vert_{L^2(1-2\eps^\gamma,1)}\\
&\leq \eps^{-1/2}\biggl(\vert y_0(1)\vert \Vert e^{-Mz}\Vert_{L^2(1-2\eps^\gamma,1)}+\vert y^{(1)}_0(1)\vert \Vert z e^{-Mz}\Vert_{L^2(1-2\eps^\gamma,1)} \biggr)\\
& \leq c (\vert y_0(1)\vert+ \vert y^{(1)}_0(1)\vert).
\end{aligned}
$$
In particular, if the initial condition and its first derivative vanish at $x=1$, then we have the following $L^2(H^1)$ estimate :
\begin{theorem}
\label{convergenceL2H1}
Let $y^{\eps}$ be the solution of \eqref{eq:2.1} and $\widetilde{P}^{\eps}$ the function defined by \eqref{def.P3epstilde}. Assume \eqref{assumption} and that $y_0(1)=y_0^{(1)}(1)=0$. Then there exists a constant $c>0$ independent of $\eps$, such that, 
\begin{equation}
\label{rateL2H1}
 \Vert (\widetilde{P}^{\eps}-y^{\eps})_x \Vert_{L^2(Q_T)} \leq c\,\eps. 
 \end{equation}
\end{theorem}
In the general case for which $y_0(1)\neq 0$ and $y_0^{(1)}(1)\neq 0$, we may achieve the same rate by making an asymptotic analysis of the corner layer of the solution $y^{\eps}$ at the point $(x,t)=(1,0)$: as mentioned earlier, this requires to introduce the auxiliary scale variable $\tau_1=t/\sqrt{\eps}$.
\end{remark}

\begin{remark}
The negative case $M<0$ (leading to surprising results for the corresponding null controllability problem, see \cite{CoronGuerrero2005,munchSIMAI2019}) exhibits a boundary layer at $x=0$ and an internal layer along the second characteristic $\{(x,t)\in Q_T, x+Mt-1=0 \}$. This case can be treated in a similar way. Actually, using the change of variable $\overline{x}=1-x$, we see that $\overline{y}^{\eps}(\overline{x},t)=y^\eps(x,t)$ solves the advection-diffusion equation 
\begin{equation*}
\left\{
\begin{aligned}
& \overline{y}^{\varepsilon}_t - \varepsilon \overline{y}^{\varepsilon}_{\overline{x}\overline{x}} +(-M) \overline{y}^{\varepsilon}_{\overline{x}}=0, & (x,t) \in Q_T,\\
& \overline{y}^{\varepsilon}(0,t)=0,  \quad \overline{y}^{\varepsilon}(1,t)=v(t), & t\in (0,T),\\
& \overline{y}^{\varepsilon}(x,0)=y_0(1-x), & x\in (0,1).
\end{aligned}
\right.
\end{equation*}
Since now $-M>0$, it suffices to adapt the analysis of the previous sections by interchanging the Dirichlet conditions. Expressions of the functions $y^k(x,t)$ are simpler since they vanishes above the first characteristic. On the other hand, the functions $Y^k(z,\tau,t)$ do not vanish anymore at $z=0$. 
\end{remark}

\section{Application: Estimate of $\Vert y^\eps(\cdot,1/M)\Vert_{L^2(0,1)}$ for $v\equiv 0$}\label{1surM}

The asymptotic analysis developed in the previous sections provides some information of the solution $y^{\eps}$ of (\ref{eq:transport}) associated to $v\equiv 0$, at $t=1/M$. For $t>1/M$, we recall the following decay property (see \cite{amirat_munch}) obtained using energy estimates. 
\begin{lemma}
Let $\alpha\in [0,1)$. The solution $y^{\eps}$ associated to $v\equiv 0$ satisfies 
\begin{equation*}
\Vert y^{\eps}(\cdot,t)\Vert_{L^2(0,1)}\leq \Vert y^{\eps}(\cdot,0)\Vert_{L^2(0,1)}e^{-\frac{M\alpha^2}{4\eps(1-\alpha)}}, \quad \forall t\geq \frac{1}{M(1-\alpha)}.
\end{equation*}
\end{lemma}
The $L^2(0,1)$-norm of the solution at any time strictly greater than $1/M$ is therefore exponentially small with respect to $\eps$. This is the effect of the transport term. At $t=1/M$, the behavior of the $L^2(0,1)$-norm with respect to $\eps$ is polynomial with a rate which depends on the derivatives of the initial condition $y_0$ at $x=0$. We have the following estimate.

\begin{proposition}\label{polynomial_decay}
Assume $v\equiv 0$. For $\eps>0$ small enough, the solution $y^\eps$ of (\ref{eq:transport}) satisfies 
\begin{equation}
\biggl\Vert y^\eps\biggl(\cdot,\frac{1}{M}\biggr)\biggr\Vert_{L^2(0,1)}\leq C\biggl(\vert y_0(0)\vert \eps^{1/4}+\vert y^{(1)}_0(0)\vert \eps^{3/4}+\vert y^{(2)}_0(0)\vert \eps^{5/4}\biggr) + \mathcal{O}(\eps^{3/2}) \label{polynomialdecay}
\end{equation}
for some constant $C>0$.
\end{proposition}
\textsc{Proof-} Using Theorem \ref{convergence}, we write that 
$$
\begin{aligned}
\Vert y^{\eps}(\cdot,1/M)\Vert_{L^2(0,1)} & \leq   \Vert \widetilde{P}^{\eps}(\cdot,1/M)\Vert_{L^2(0,1)} + \Vert (y^\eps-\widetilde{P}^{\eps})(\cdot,1/M)\Vert_{L^2(0,1)}\\
& \leq   \Vert \widetilde{P}^{\eps}(\cdot,1/M)\Vert_{L^2(0,1)} + \mathcal{O}(\eps^{3/2}).
\end{aligned}
$$ 

We then estimate the norm of $\widetilde{P}^{\eps}$ at $t=1/M$. Remark that the variable $\tau=\frac{1/M-t}{\sqrt{\eps}}$ vanishes for $t=1/M$. From (\ref{def.P3epstilde}), we compute that 
\begin{equation*}
\widetilde{P}^{\eps}\biggl(x,\frac{1}{M}\biggr)=  \sum_{k=0}^3 \eps^{\frac{k}{2}}P^{k/2}_\eps\biggl(x,\frac{1}{M}\biggr)=m_0(x)+\eps^{1/2}m_1(x)+\eps m_2(x) +\eps^{3/2}m_3(x)
\end{equation*}
where for $k=0,\cdots,3$, the function $m_k$ is factor of $y_0^{(k)}(0)$ and is defined as follows (recall that $z=(1-x)/\eps$) : 
\begin{equation*}
\begin{aligned}
& m_0(x)= W_\eps^0\biggl(\frac{x-1}{\sqrt{\eps}},\frac{1}{M}\biggr)  - \biggl[W_\eps^0\biggl(0,\frac{1}{M}\biggr)+\eps^{1/2}z W_{\eps,w}^0\biggl(0,\frac{1}{M}\biggr)\\
& \hspace*{6cm}+\eps \frac{z^2}{2} W_{\eps,ww}^0\biggl(0,\frac{1}{M}\biggr)+\eps^{3/2}\frac{z^3}{6} W_{\eps,www}^0\biggl(0,\frac{1}{M}\biggr)\biggr]e^{-Mz},\\
& m_1(x)= W^{1/2}\biggl(\frac{x-1}{\sqrt{\eps}},\frac{1}{M}\biggr) - \biggl[W^{1/2}\biggl(0,\frac{1}{M}\biggr)+\eps^{1/2}z W_w^{1/2}\biggl(0,\frac{1}{M}\biggr)+\eps \frac{z^2}{2} W_{ww}^{1/2}\biggl(0,\frac{1}{M}\biggr)\biggr]e^{-Mz},\\
& m_2(x)= W^{1}\biggl(\frac{x-1}{\sqrt{\eps}},\frac{1}{M}\biggr) - \biggl[W^{1}\biggl(0,\frac{1}{M}\biggr)+\eps^{1/2}z W_w^{1}\biggl(0,\frac{1}{M}\biggr)\biggr]e^{-Mz},\\
& m_3(x)= W^{3/2}\biggl(\frac{x-1}{\sqrt{\eps}},\frac{1}{M}\biggr) - W^{3/2}\biggl(0,\frac{1}{M}\biggr)e^{-Mz}.
\end{aligned}
\end{equation*}
$\bullet$ We first estimate $\Vert m_0\Vert_{L^2(0,1)}$.
Using that the asymptotic behavior (\ref{behavior_erf}) of the error function, we obtain
\begin{equation*}
\begin{aligned}
& W_\eps^0\biggl(\frac{x-1}{\sqrt{\eps}},\frac{1}{M}\biggr)=\frac{y_0(0)}{2}\biggl(1+erf\biggl(\frac{\sqrt{M}(x-1)}{2\sqrt{\eps}}\biggr)-e^{\frac{Mx}{\eps}}erfc\biggl(\frac{\sqrt{M}}{\sqrt{\eps}}\frac{x+1}{2}\biggr)\biggr),\\
& W_\eps^0\biggl(0,\frac{1}{M}\biggr)=\frac{y_0(0)}{2}\biggl(1-e^{\frac{M}{\eps}}erfc\biggl(\frac{\sqrt{M}}{\sqrt{\eps}}\biggr)\biggr)=\frac{y_0(0)}{2}+\mathcal{O}(\sqrt{\eps}),\\
& W_{\eps,w}^0\biggl(0,\frac{1}{M}\biggr) = y_0(0)\biggl[ \frac{\sqrt{M}}{\sqrt{\pi}}-\frac{M}{2\sqrt{\eps}}e^{\frac{M}{\eps}}erfc\biggl(\frac{\sqrt{M}}{\sqrt{\eps}}\biggr)\biggr]= y_0(0)\frac{\sqrt{M}}{2\sqrt{\pi}}+\mathcal{O}(\eps),\\
& W_{\eps,ww}^0\biggl(0,\frac{1}{M}\biggr)=\frac{y_0(0)}{2}\biggl(\frac{M^{3/2}}{\sqrt{\eps\pi}}-\frac{M^2}{\eps}e^{\frac{M}{\eps}}erfc\biggl(\frac{\sqrt{M}}{\sqrt{\eps}}\biggr)\biggr)=y_0(0)\frac{\sqrt{M}}{4\sqrt{\pi}}\eps^{1/2} +\mathcal{O}(\eps^{3/2}),\\
& W_{\eps,www}^0\biggl(0,\frac{1}{M}\biggr)=\frac{y_0(0)}{2}\biggl(-\frac{M^{3/2}}{\sqrt{\pi}}+\frac{M^{5/2}}{\eps\sqrt{\pi}}-\frac{M^3}{\eps^{3/2}}e^{\frac{M}{\eps}}erfc\biggl(\frac{\sqrt{M}}{\sqrt{\eps}}\biggr)\biggr)=-y_0(0)\frac{M^{3/2}}{4\sqrt{\pi}}+\mathcal{O}(\eps).
\end{aligned}
\end{equation*}
Recalling that that $z(x)=(1-x)/\eps$ and that $\Vert z(x)^n e^{-Mz(x)}\Vert_{L^2(0,1)}=\mathcal{O}(\sqrt{\eps})$ for all $n\in \mathbb{N}$, we deduce that the $L^2(0,1)$-norm of 
$$
w_3(x):=- \biggl[W_\eps^0\biggl(0,\frac{1}{M}\biggr)+\eps^{1/2}z W_{\eps,w}^0\biggl(0,\frac{1}{M}\biggr)+\eps \frac{z^2}{2} W_{\eps,ww}^0\biggl(0,\frac{1}{M}\biggr)+\eps^{3/2}\frac{z^3}{6} W_{\eps,www}^0\biggl(0,\frac{1}{M}\biggr)\biggr]e^{-Mz}
$$
is 
$\Vert w_3\Vert_{L^2(0,1)}=\vert y_0(0)\vert \mathcal{O}(\eps^{1/2})$. We then write that 
\begin{equation*}
\begin{aligned}
m_0(x)=&\frac{y_0(0)}{2}\biggl[\underbrace{1+erf\biggl(\frac{\sqrt{M}(x-1)}{2\sqrt{\eps}}\biggr)}_{:=w_1(x)}\underbrace{-e^{\frac{Mx}{\eps}}erfc\biggl(\frac{\sqrt{M}}{\sqrt{\eps}}\frac{x+1}{2}\biggr)}_{:=w_2(x)}\biggr]+w_3(x).
\end{aligned}
\end{equation*}
Using (\ref{estimate_chu}) with $y=\sqrt{M}(1-x)/ (2\sqrt{\eps})\geq 0$, we obtain 
$$
\frac{1}{2}e^{-\frac{M(1-x)^2}{\pi\eps}}\leq w_1(x)\leq e^{-\frac{M(1-x)^2}{4\eps}}, \quad\forall x\in [0,1]
$$
leading to $\Vert w_1\Vert_{L^2(0,1)}=\mathcal{O}(\eps^{1/4})$. Moreover, using the asymptotic behavior of the error function $erf$, we have  
$$
w_2(x)=-e^{-\frac{M}{\eps}\frac{(1-x)^2}{4}}\biggl(\frac{\sqrt{\eps}}{\sqrt{M\pi}}\frac{2}{x+1}+\mathcal{O}(\eps^{3/2})\biggr), \quad\forall x\in [0,1].
$$
We then deduce that $\Vert w_2\Vert_{L^2(0,1)}=\mathcal{O}(\eps^{3/4})$. We also check that $\int_0^1 w_1(x)w_2(x)dx=\mathcal{O}(\eps)$, $\int_0^1 w_2(x)w_3(x)dx=\mathcal{O}(\eps^{3/2})$ and $\int_0^1 w_1(x)w_3(x)dx=\mathcal{O}(\eps)$. This allows to conclude that 
$$
\Vert m_0\Vert_{L^2(0,1)}=\vert y_0(0)\vert \mathcal{O}(\eps^{1/4}).
$$
$\bullet$ We now estimate $\eps^{1/2}\Vert m_1\Vert_{L^2(0,1)}$. We have 
\begin{equation*}
\begin{aligned}
& W^{1/2}\biggl(\frac{x-1}{\sqrt{\eps}},\frac{1}{M}\biggr) =\frac{y_0^{(1)}(0)}{2} \biggl[ w_1(x)\frac{(x-1)}{\sqrt{\eps}}  + \frac{2}{\sqrt{M\pi}} e^{-\frac{M (1-x)^2}{4\eps}}  \biggr],\\
& W^{1/2}\biggl(0,\frac{1}{M}\biggr) =\frac{y_0^{(1)}(0)}{\sqrt{M\pi}},\quad  W_w^{1/2}\biggl(0,\frac{1}{M}\biggr) =\frac{y_0^{(1)}(0)}{2},\quad W_{ww}^{1/2}\biggl(0,\frac{1}{M}\biggr) =\frac{\sqrt{M}}{\sqrt{\pi}}\frac{y_0^{(1)}(0)}{2}.
\end{aligned}
\end{equation*}
Using the estimate above for $w_1$, we check that  $\Vert W^{1/2}\biggl(\frac{x-1}{\sqrt{\eps}},\frac{1}{M}\biggr)\Vert_{L^2(0,1)}=\mathcal{O}(\eps^{1/4})$. Its follows that $\Vert m_1\Vert_{L^2(0,1)}=\vert y^{(1)}_0(0)\vert  \mathcal{O}(\eps^{1/4})$ and 
$$
\eps^{1/2}\Vert m_1\Vert_{L^2(0,1)}=\vert y^{(1)}_0(0)\vert  \mathcal{O}(\eps^{3/4}).
$$
$\bullet$ We now estimate $\eps\Vert m_2\Vert_{L^2(0,1)}$.
\begin{equation*}
\left\{
\begin{aligned}
& W^1\biggl(\frac{x-1}{\sqrt{\eps}},\frac{1}{M}\biggr)= \frac{y^{(2)}_0(0)}{2}\biggl[\biggl(\frac{(1-x)^2}{2\eps}+\frac{1}{M}\biggr)w_1(x) + \frac{x-1}{\sqrt{M\pi}\sqrt{\eps}} e^{-\frac{M(1-x)^2}{4\eps}}\biggr],\\
& W^1\biggl(0,\frac{1}{M}\biggr)= \frac{y^{(2)}_0(0)}{2M},\quad W_w^1\biggl(0,\frac{1}{M}\biggr)= \frac{y^{(2)}_0(0)}{\sqrt{M\pi}}.
\end{aligned}
\right.
\end{equation*} 
Using the estimate above for $w_1$, we again that  $\Vert W^1\big(\frac{x-1}{\sqrt{\eps}},\frac{1}{M}\big)\Vert_{L^2(0,1)}=\mathcal{O}(\eps^{1/4})$. Its follows that $\Vert m_2\Vert_{L^2(0,1)}=\vert y^{(1)}_0(0)\vert  \mathcal{O}(\eps^{1/4})$ and 
$$
\eps\Vert m_2\Vert_{L^2(0,1)}=\vert y^{(2)}_0(0)\vert  \mathcal{O}(\eps^{5/4}).
$$
Similarly, we compute $\eps^{3/2}\Vert m_3\Vert_{L^2(0,1)}=\vert y^{(3)}_0(0)\vert  \mathcal{O}(\eps^{7/4})$. For $\eps>0$ small enough, this term is absorbed by $\Vert (y^{\eps}-\widetilde{P}^{\eps})(\cdot,1/M)\Vert_{L^2(0,1)}$.  Then, writing that $\Vert \widetilde{P}^\eps(\cdot,1/M)\Vert_{L^2(0,1)}=\Vert \sum_{k=0}^3 \eps^{k/2}m_k\Vert_{L^2(0,1)}\leq \sum_{k=0}^3 \Vert\eps^{k/2}m_k\Vert_{L^2(0,1)}$, we obtain the result. $\hfill\Box$

Under additional regularity assumptions, estimate (\ref{polynomialdecay}) suggests that the norm $\Vert y^{\eps}(\cdot,1/M)\Vert_{L^2(0,1)}$ decays exponentially with respect to $\eps$
if the derivative $y_0^{(j)}(0)$ vanishes for all $j\in \mathbb{N}$.

As an illustration, we consider the simple case $v\equiv 0$ and $y_0\equiv 1$ for which 
$$
\left\{
\begin{aligned}
&\widetilde{P}^\eps(x,t)= W_\eps^0(w,t)-\biggl(W^0_{\eps}(M\tau,t)+\eps^{1/2}z W^0_{\eps,w}(M\tau,t)+\\
&\hspace*{5cm}\eps \frac{z^2}{2} W^0_{\eps,ww}(M\tau,t)+\eps^{3/2} \frac{z^3}{6} W^0_{\eps,www}(M\tau,t)\biggr)e^{-Mz},\\
& w=\frac{x-Mt}{\sqrt{\eps}}, \quad M\tau= \frac{1-Mt}{\sqrt{\eps}}, \quad z=\frac{1-x}{\eps}.
\end{aligned}
\right.
$$
Figure \ref{widetildeP} depicts the function $\widetilde{P}^{\eps}(x,t)$ over $x\in (0,1)$ for $t=1/(2M)$ and $t=1/M$. We take $M=1$ and $\eps\in \{10^{-2},10^{-3}\}$. As $\eps$ goes to zero, the function 
$\widetilde{P}^{\eps}(x,1/(2M))$ displays a transition from $0$ to $1$ at the point $x=1/2$ and a faster transition from $1$ to $0$ at the point $x=1^-$. For $t=1/M$, these two transitions, from $0$ to $1$ and from $1$ to $0$ occur simultaneously in the neighborhood of $x=1$.  Figure \ref{plot3d} depicts the approximation $\widetilde{P}^\eps$, for $\eps=10^{-2}$, $M=1$ in the domain $(0,1)\times (0,1.2/M)$ and clearly show the simultaneous occurence of both an internal and boundary layers.

\begin{figure}[!http]
\begin{center}
\includegraphics[scale=0.34]{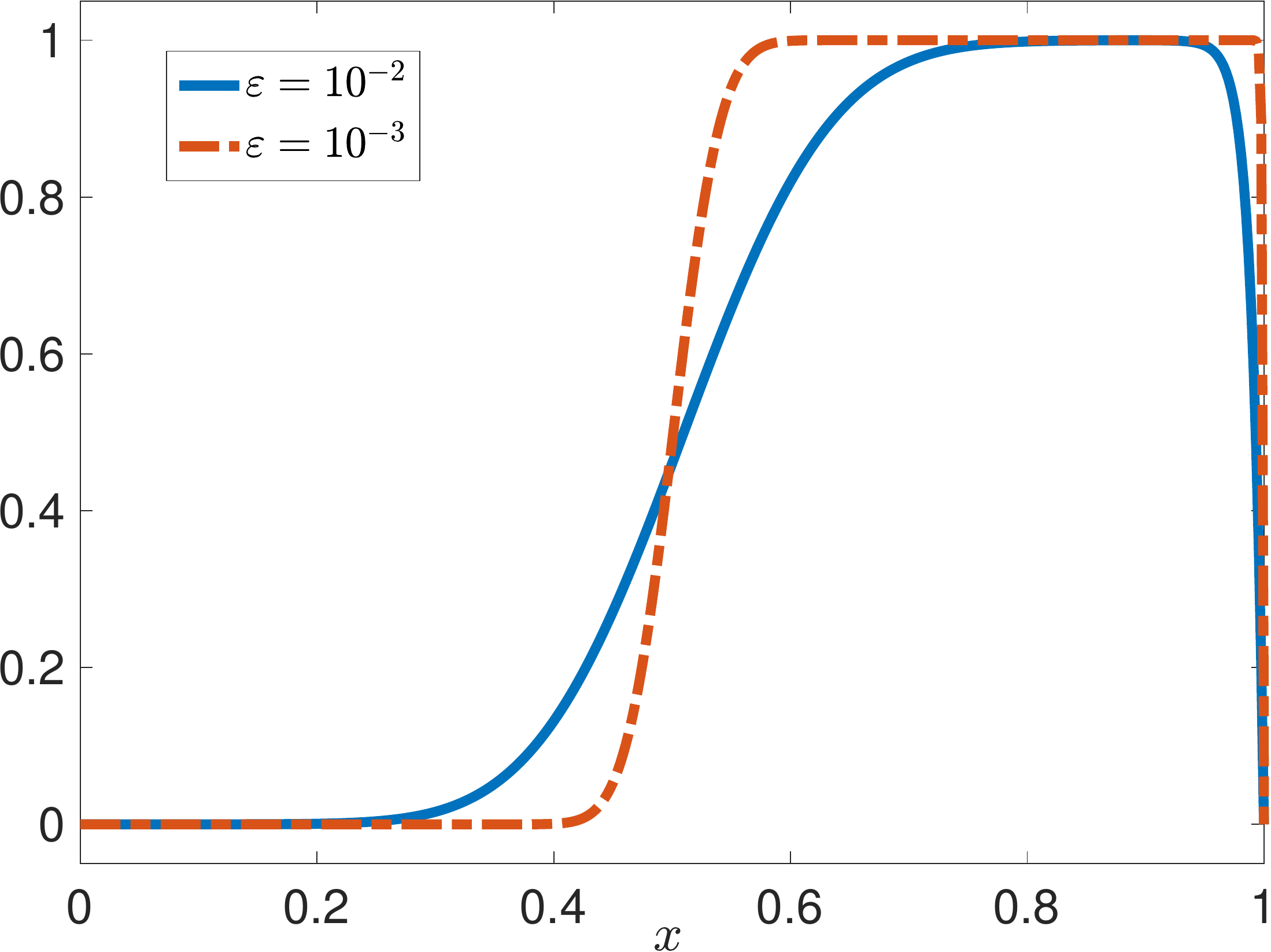}\hspace*{0.7cm}
\includegraphics[scale=0.34]{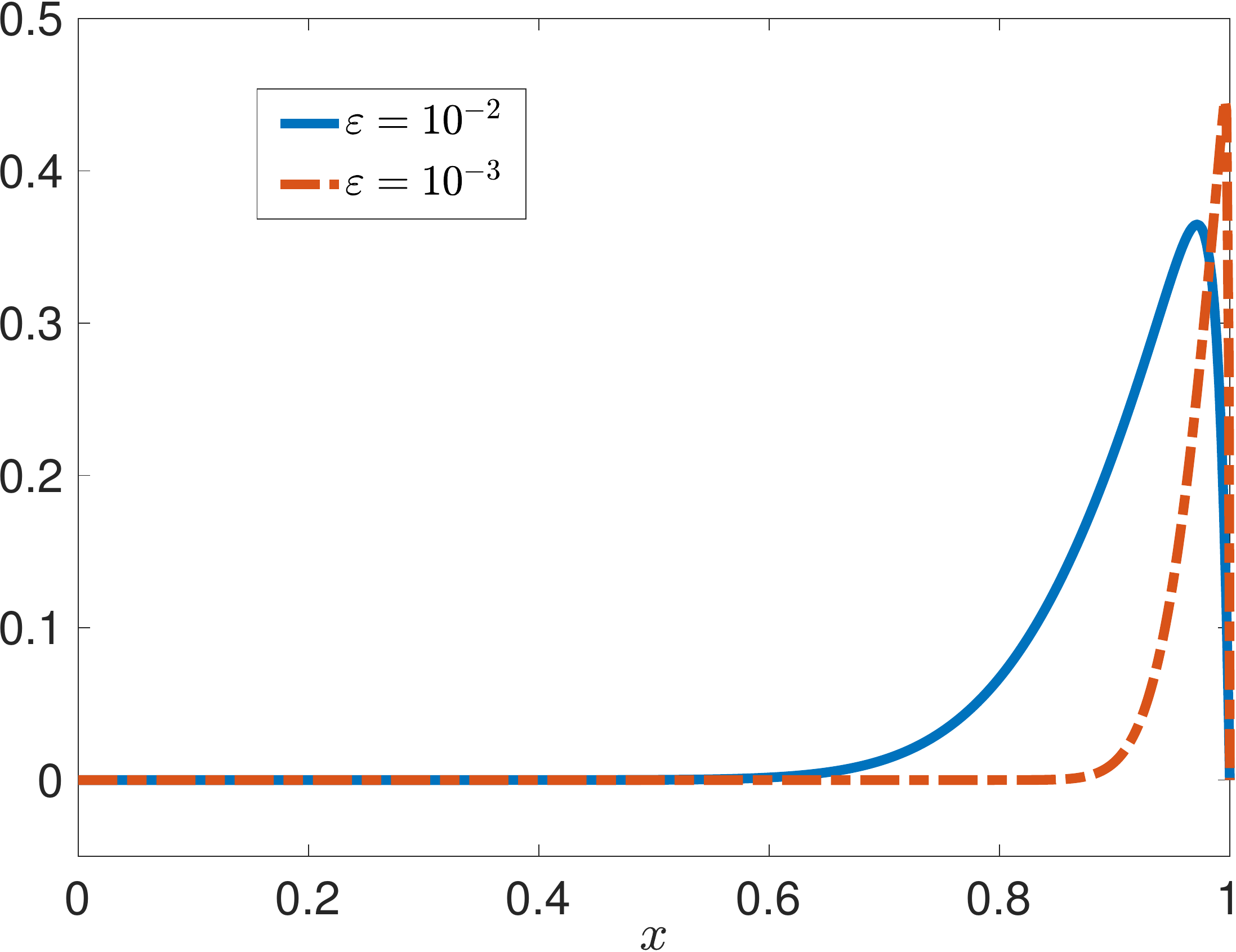}
\caption{$\widetilde{P}^{\eps}(x,t)$ w.r.t. $x\in (0,1)$ for $t=1/2/M$ (left) and $t=1/M$ (right); $M=1$, $\eps\in \{10^{-2},10^{-3}\}$.}\label{widetildeP}
\end{center}
\end{figure}

\begin{figure}[!http]
\begin{center}
\includegraphics[scale=0.5]{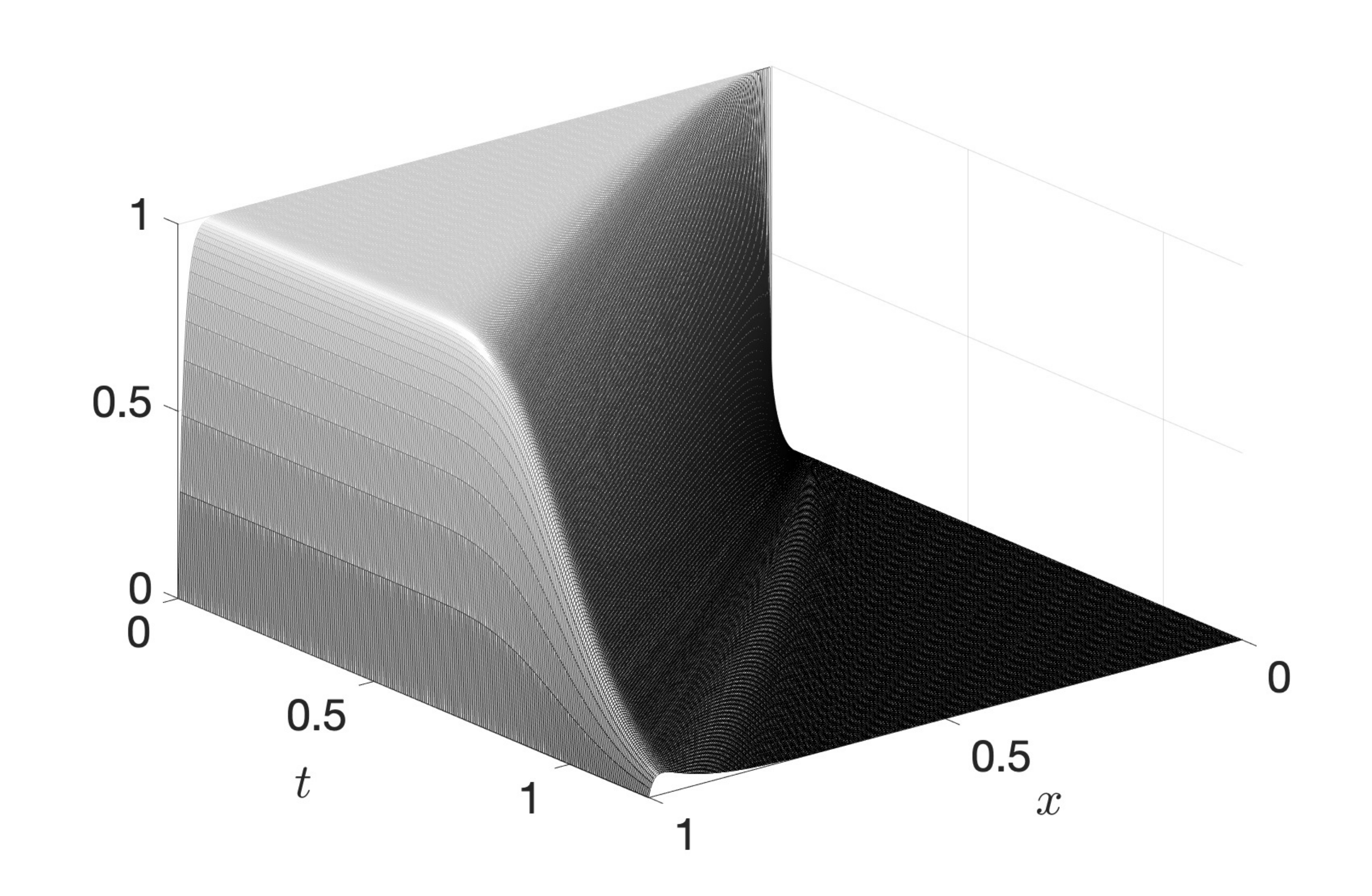}
\caption{$\widetilde{P}^{\eps}(x,t)$  in $(0,1)\times (0,1.2/M)$; $M=1$, $\eps=10^{-2}$; $v\equiv 0$, $y_0\equiv 1$.}\label{plot3d}
\end{center}
\end{figure}

\section{Concluding remarks and perspectives}\label{conclusions}

We have shown that the method of matched asymptotic expansion with appropriate scaling allows to approximate solutions of a singular boundary value 
problem involving interacting internal and boundary layers of distinct sizes. The approximation is a linear combination of three expansions, describing the three
behavior of the solution $y^\eps$ as $\eps$ goes to zero. These expansions are locally matched and lead to a $C^1(Q_T)$-approximation.  Using standard energy estimates, and assuming regularity on the data, precisely, $v\in C^4([0,T])$ and $y_0\in C^4([0,1])$, this approximation, denoted by $\widetilde{P}^\eps$, notably fulfills the property $\Vert y^\eps-\widetilde{P}^\eps\Vert_{L^\infty(0,T; L^2(0,1))}\leq C \eps^{3/2}$. With more regularity, we can actually achieve an arbitrarily large rate. This requires however a large amount of calculus. 

Similarly, the method can be extanded to the case $\Omega\in \mathbb{R}^2$ and the case of non constant coefficient (we refer to \cite{plaschko1990,Shih2001} where the equation $y^{\eps}_t -\eps y^\eps_{xx}+M(t)y^{\eps}_x=0$, for $x\in \mathbb{R}^+$ is analyzed). It would also be interesting to consider the case of nonlinear equations like the Burger's equation $y^{\eps}_t-\eps y^{\eps}_{xx}+y^{\eps} y^{\eps}_x=0$ introduced to model turbulence. The asymptotic analysis of this equation posed for $x\in \mathbb{R}$ is mentioned in \cite[Section 4.3.1]{ColeBook81} (we also refer to \cite{Marbach2014}). 

Eventually, we mention that our analysis assumes that the initial condition is independent of the parameter $\eps$. In particular, our analysis does not apply for the initial condition $y_0^{\eps}(x)=K_\eps e^{-\frac{Mx}{2\eps}}\sin(\pi x)$, with $K_{\eps}=\mathcal{O}(\eps^{-3/2})$ so that $\Vert y_0^\eps\Vert_{L^2(0,1)}=1$ exhibited in \cite{munchSIMAI2019}. As $\eps$ goes to 0, this initial condition gets concentrated at $x=0$ and is suspected to maximize the corresponding cost of null control for (\ref{eq:transport}), discussed in \cite{CoronGuerrero2005}. However, by introducing the new function $z^{\eps}(x,t):=K_{\eps}^{-1}e^{\frac{Mx}{2\eps}}y^\eps(x,t)$, we check that $z^{\eps}$ solves (taking $v\equiv 0$) the boundary value problem 
\begin{equation*}
\left\{
\begin{aligned}
& z^{\varepsilon}_t - \varepsilon z^{\varepsilon}_{xx} + 2M z^{\varepsilon}_x - \frac{M^2}{4\eps}z^\eps=0, & (x,t) \in Q_T,\\
& z^{\varepsilon}(0,t)=0,  \quad z^{\varepsilon}(1,t)=0, & t\in (0,T),\\
& z^{\varepsilon}(x,0)=sin(\pi x), & x\in (0,1),
\end{aligned}
\right.
\end{equation*}
on which we may apply our asymptotic analysis to get an approximation of $z^{\eps}$ and then obtain notably the order of magnitude of the norm $\Vert y^{\eps}(\cdot,1/M)\Vert_{L^2(0,1)}$.

\bibliographystyle{siam}

\bibliography{biblio_transport_Asymptotic.bib}

\end{document}